\documentclass[12pt]{amsart}

\usepackage{amsmath,amscd}
\usepackage{amssymb}
\usepackage{amsthm}
\usepackage{oldgerm}
\usepackage{multicol}

\usepackage[pdftex,colorlinks,urlcolor=blue,pdfstartview=FitH]{hyperref}

\newcommand{\val}{{\sf val}}

\newcommand{\ZZ}{\mathbb{Z}}
\newcommand{\R}{\mathbb{R}}

\newcommand{\As}{\mathcal{A}s}

\newcommand{\mQ}{\ensuremath{\mathbb{Q}}}

\newcommand{\mG}{\ensuremath{\mathcal{G}}}

\newcommand{\End}{{\sf End}}
\newcommand{\Hom}{{\sf Hom}}
\newcommand{\Tw}{{\sf Tw}}

\newcommand{\Zm}{\ensuremath{\mathbb{Z}}}

\newcommand{\mM}{\ensuremath{\mathcal{M}}}

\newcommand{\mA}{\ensuremath{\mathcal{A}}}
\newcommand{\mO}{\ensuremath{\mathcal{O}}}

\newcommand{\LL}{\mathcal{L}}

\newcommand{\id}{{\sf id}}
\newcommand{\BB}{{\mathbb{B}}}
\newcommand{\Codiff}{{\sf Codiff}}

\newtheorem{thm}{Theorem}[section]
\newtheorem{lem}[thm]{Lemma}

\newtheorem{cor}[thm]{Corollary}
\newtheorem{prop}[thm]{Proposition}
\newtheorem{rmk}[thm]{Remark}
\newtheorem{constr}[thm]{Construction}

\theoremstyle{definition}
\newtheorem{defn}[thm]{Definition}

\def\proof {\noindent{\sc{Proof. }}}
\def\qed {\mbox{}\hfill {\small \fbox{}} \\}
\def\lto{\longrightarrow}

\def\leq{\leqslant}
\def\geq{\geqslant}
\def\Ai{A$_\infty$}
\newcommand{\Aii}{\ensuremath{\mA_\infty}}

\newcommand{\m}{\mathfrak{m}}
\def\Ai{$A_\infty$}

\title{Tensor Product of Cyclic A$_\infty$-Algebras and their Kontsevich Classes}

\author {LINO AMORIM \textsuperscript{1}}
\thanks{\textsuperscript{1}Mathematical Institute, University of Oxford, Andrew Wiles Building, Radcliffe Observatory Quarter, Woodstock Road, Oxford, OX2 6GG, England. {\tt Email: camposamorim@maths.ox.ac.uk}}
\author{JUNWU TU \textsuperscript{2}}
\thanks{\textsuperscript{2}Department of Mathematics, University of Missouri, Columbia, 65211, USA.  {\tt Email: tuju@missouri.edu}}
\date{}

\begin{document}

\begin{abstract}

Given two cyclic  $A_\infty$-algebras $A$ and $B$, in this paper we prove that there exists a cyclic $A_\infty$-algebra structure on their tensor product $A\otimes B$ which is unique up to a cyclic $A_\infty$-quasi-isomorphism. Furthermore, the Kontsevich class of  $A\otimes B$ is equal to the cup product of the Kontsevich classes of $A$ and $B$ on the moduli space of curves.

\end{abstract}

\maketitle

\section{Introduction}

\subsection{Tensor products} Let $A$ and $B$ be two associative algebras with product maps $\m_2^A$ and $\m_2^B$. Then there is a natural associative algebra structure on $A\otimes B$ defined by
\[ \m_2^\otimes (a_1\otimes b_1, a_2\otimes b_2):=\m_2^A(a_1,a_2)\otimes m_2^B(b_1,b_2).\]
This tensor product construction is much less obvious if one replaces associative algebras by their homotopy version: {\sl $A_\infty$-algebras}, introduced by Stasheff~\cite{Sta}. Explicitly, this structure consists of a  graded vector space $A$ and a family of of multi-linear maps of degree $k-2$
$$\m_k: A^{\otimes k} \rightarrow A$$
for each $k\geq 1$, which satisfy a certain homotopy version of associativity. In this case, the naive definition of a tensor product structure
\[ \m_k^\otimes:=\m_k^A\otimes \m_k^B\]
is clearly false, simply by degree considerations.

In~\cite{SanUmb}, Saneblidze and Umble constructed a tensor product structure on $A\otimes B$ whose structure maps $m_k^\otimes$ are given by a beautiful formula involving various compositions of $\m_j$ $(j\leq k)$ in each tensor component. Later in~\cite{MarShn}, Markl and Shnider interpreted the Saneblidze-Umble's formula as giving a diagonal map on the  associahedra (the cellular complex which governs $A_\infty$-algebras). Their construction used cubical decompositions of the associahedra. A similar approach, using simplicial decompositions of the associahedra instead, was carried out by Loday~\cite{Lod}.

In this paper, we are interested in the tensor product of $A_\infty$-algebras with an additional structure: a cyclic inner product. Namely, a {\em cyclic} $A_\infty$-algebra is an $A_\infty$-algebra $A$, together with a non-degenerate inner product $\langle,\rangle$ such that the expression
\[ \langle \m_k(a_0,\ldots,a_{k-1}),a_k\rangle\]
is invariant under cyclic permutations of the indices $0,1,\ldots,k$ (up to a sign). 

In the case when the only nonzero structure map is $\m_2$, this is equivalent to the notion of a (noncommutative) Frobenius algebra. Given two Frobenius algebras $A$ and $B$, the natural tensor product algebra structure $\m_2^\otimes$ on $A\otimes B$ is again Frobenius with respect to the natural inner product defined by
\begin{equation}~\label{pairing-eq}
\langle a_1\otimes b_1, a_2\otimes b_2\rangle_{A\otimes B}:= (-1)^{|b_1||a_2|}\langle a_1, a_2\rangle_A \langle b_1, b_2\rangle_B .\end{equation}
Thus a natural question is whether Saneblidze-Umble's formula $\m_k^\otimes$ on the tensor product is cyclic with respect to the pairing above.
Unfortunately, this is not the case. As was observed in~\cite{TraUmb}, the third product $\m_3^\otimes$ of Saneblidze-Umble is already not cyclic. 

In \cite{Cho}, Cho introduced the notion of {\it strong homotopy inner product}, which is a cyclic \Ai-algebra up to homotopy. It consists of an infinity inner product as defined by Tradler (see \cite{TraUmb}, for a definition) satisfying some additional properties, in particular being {\it closed}. He also shows that an \Ai-algebra with a strong homotopy inner product, is quasi-isomorphic to a cyclic \Ai-algebra. In ~\cite{TraUmb}, the authors show that the tensor product of cyclic \Ai-algebras admits an infinity inner product, but did not check if it is closed. Therefore we cannot use Cho's result to conclude the tensor product is cyclic up to quasi-isomorphism. 

In this paper, we take a different approach to this problem. Namely, we use a different definition of tensor product and we prove the following

\medskip
\begin{thm}~\label{main-thm}
Let $A$ and $B$ be two cyclic $A_\infty$-algebras. Then there exists a cyclic $A_\infty$-algebra structure on $A\otimes B$ with respect to the inner product defined by Formula~(\ref{pairing-eq}). Furthermore, this structure is unique up to cyclic homotopy.
\end{thm}

The theorem is proved by showing that the cyclic differential graded operad $\Aii$ admits a {\em cyclic} diagonal morphism
\[ \Delta: \Aii \rightarrow \Aii\otimes \Aii\]
which is unique up to homotopy of morphisms between cyclic operads. Our proof is to inductively construct $\Delta$ by using the simple fact that the associahedra are contractible. Unlike the case of tensor products without cyclic structure, we do not obtain a general formula for $\m_k^\otimes$, but we give in Section~\ref{formula-sub} a description of $\m_k^\otimes$ for $k\leq4$.

\subsection{Kontsevich classes} It is well-known that the associative graph homology complex $\mG_*$ with $\mQ$ coefficients computes the rational homology of $\coprod_{g,n} \mM_{g,n}$, the union of the moduli spaces of smooth Riemann surfaces of genus $g$ and $n$ unlabeled boundary components . In~\cite{Kon}, Kontsevich associated a cocycle $c_A$ on the  complex $\mG_*$, to any given cyclic $A_\infty$-algebra, with finite dimensional cohomology, whose inner product is even and symmetric.  Thus the cocycle $c_A$ also defines a cohomology class $[c_A]\in H^*(\coprod \mM_{g,n},K)$ if the algebra $A$ is defined over $K$, a field of characteristic zero. 

In the case when the inner product is odd, but still symmetric, one can define a class $[c_A]\in H^*(\coprod \mM_{g,n}, \LL)$ where $\LL$ is a natural local system on $\coprod \mM_{g,n}$ whose fiber over a (connected) Riemann surface $\Sigma$ of type $(g,n)$ is
\[ \LL|_{\Sigma}:=\det H^1( \Sigma, K),\]
the determinant of the first cohomology group of $\Sigma$.
Since $H^1(\Sigma,K)\cong H^1(\Sigma,\mathbb{Z})\otimes K$, it is clear that $\LL^{\otimes 2}$ is trivial. 
In summary, for a cyclic $A_\infty$ algebra $A$ whose inner product is of parity $\epsilon\in \mathbb{Z}/2\mathbb{Z}$, there is an associated class $c_A\in H^*(\coprod \mM_{g,n}, \LL^{\otimes \epsilon})$. Note that since $\LL^{\otimes 2}$ is trivial, the notation $\LL^{\otimes \epsilon}$ causes no confusion. 

We can now state our second main

\medskip
\begin{thm}~\label{main2-thm}
Let $A$ and $B$ be two finite dimensional, cyclic $A_\infty$-algebras over a field $K$ of characteristic zero. Assume that the inner products on $A$ and $B$ are both symmetric and of parity $\epsilon_1,\epsilon_2\in \mathbb{Z}/2\mathbb{Z}$ respectively. Then we have
\[ [c_{A\otimes B}] = [c_A]\cup [c_B]\]
where $\cup: H^*(\coprod \mM_{g,n},\LL^{\otimes \epsilon_1})\otimes H^*(\coprod \mM_{g,n},\LL^{\otimes \epsilon_2}) \rightarrow H^*(\coprod \mM_{g,n},\LL^{\otimes\epsilon_1+\epsilon_2})$ is the cup product map on cohomology with local coefficients.
\end{thm}

The main ingredient in the proof of this result is the construction of a cellular ``diagonal"
$$\delta: \mG_*(\LL ^{\otimes \epsilon_1+\epsilon_2})\rightarrow \mG_*(\LL^{\otimes \epsilon_1})\otimes\mG_*(\LL^{\otimes \epsilon_2})$$
on the graph homology complex (twisted by appropriate local systems). For this, we need to fix a cyclic diagonal $\Delta$ for $\Aii$. Then we define $\delta$ applied to an oriented graph $(\Gamma,\sigma)$ to be a kind of direct product of $\Delta$ applied to every vertex of $\Gamma$. Then we show that, if one uses $\Delta$ to define the cyclic $A_\infty$-structure on $A\otimes B$, there is a chain level equality $c_{A\otimes B} = (c_A\otimes c_B)\circ \delta$,
which implies the theorem.

\subsection{Contents} The paper is organized as follows. In Section~\ref{basic-sec}, we review basic definitions about $A_\infty$-algebras and the operad $\Aii$. In Section~\ref{cyclic-sec}, we show the existence and uniqueness of a cyclic diagonal of $\Aii$, and use it to prove Theorem~\ref{main-thm}. Section~\ref{graph-sec} contains the construction of the diagonal $\delta:\mG_*\rightarrow \mG_*\otimes\mG_*$, which is then used to prove Theorem~\ref{main2-thm} in the even case. In Section~\ref{sec:twisted} we consider the case of odd \Ai-algebras.

\subsection{Conventions} Throughout the paper we work over the field $\mQ$. This is just for simplicity, all the results work over any field of characteristic $0$. For multi-linear operations with graded vector spaces, the Koszul sign convention is assumed. We work with homological complexes where the differentials are of degree $-1$. In this paper, an inner product on a graded finite dimensional vector space is always assumed to be {\sl symmetric} and either even or odd. More precisely, we require that, given homogeneous elements $a$ and $b$
$$\langle  a,b\rangle=(-1)^{\vert a\mid\mid b\mid}\langle  b,a\rangle ,$$
$$\langle  a,b\rangle \neq 0 \Rightarrow \vert a \vert \equiv \vert b \vert + \epsilon \pmod 2$$
where $\epsilon=0$ in the even case and $\epsilon=1$ in the odd case.

\subsection{Acknowledgments} The paper started with a question of Kevin Costello concerning the Kontsevich classes of the Fermat quintic. We warmly thank him for valuable discussions. The second named author would like to thank Andrei C\u ald\u araru for introducing him to the subject of $A_\infty$-algebras and graph complexes. 

During the preparation of this paper the first author was supported by EPSRC grant EP/J016950/1.

\section{The $A_\infty$-operad}~\label{basic-sec}

In this section, we recall basic definitions and facts about the operad $\Aii$. 
The book \cite{operads} is a good reference for this material.

\subsection{The $A_\infty$-operad.}

Stasheff introduced in \cite{Sta}, for each $n\geq 2$, a cellular complex $K(n)$ of dimension $n-2$. They are know as the associahedra or the \emph{Stasheff polytopes}. The  \Ai-operad is a non-symmetric, non-unital operad on the category of differential graded vector spaces, defined by setting
$$\Aii(n):= C_*(K(n),\mQ),\ n\geq 2,$$
where $C_*(K(n),\mQ)$ is the cellular chain complex of the $n$-th Stasheff polytope $K(n)$. Cells in $K(n)$ of dimension $k$ are in one-to-one correspondence with planar, rooted trees with $n$ leaves and $n-k-2$ internal edges.

To describe the differential on $\Aii(n)$ in terms of planar trees, one needs to define appropriate orientations on them. Following~\cite{MarShn}, we define an orientation of a tree $T$ with internal edges $e_1,\ldots,e_k$, to be an orientation of the vector space spanned by $\{e_1,\ldots,e_k\}$. With this notation we define $\Aii(n)$ to be the vector space generated by pairs $(T,\omega)$ where $T$ is a planar rooted tree with $n$ leaves and $\omega$ is an orientation of $T$, together with the relation $(T,-\omega)=-(T,\omega)$. The operadic composition
$$\circ_i: C_k(K(n_1))\otimes C_l(K(n_2))\lto C_{k+l}(K(n_1+n_2-1))$$
is defined by
\begin{align}\label{composition}
(U,\sigma)\circ_i(V,\tau)=(-1)^{i(n_2+1)+n_1l}(U\circ_iV,\sigma\wedge\tau\wedge e)
\end{align}
where $U\circ_iV$ is the tree obtained by grafting the root of $V$ to the $i$-th leaf of $U$ and $e$ is the new internal edge created by grafting.

The differential
$$\partial: C_k(K(n))\lto C_{k-1}(K(n))$$
is defined by
$$\partial(U,\sigma)=\sum_{\{V\mid V/e= U\}}(V, e\wedge \sigma)$$
with the sum taken over all trees $V$ with an edge $e$ such that when we collapse $e$ on $V$ we obtain $U$.

It is a well-known fact that, as an operad of graded vector spaces, $\Aii$ is freely generated by the corollas $c_k\in\Aii(k)$, the trees with only one internal vertex.

\subsection{Cyclic structures.}

The notion of a \emph{cyclic structure on an operad} $\mO$ was introduced in \cite{GetKap}. In the context of non-symmetric operads, this structure reduces to an action of the cyclic group $\Zm/(n+1)\Zm$ on $\mO(n)$ satisfying certain compatibility conditions. We refer to \cite{operads} or \cite{GetKap} for a precise definition. It was proved in \cite{GetKap} that $\Aii$ is cyclic. Denote by $r$ a generator of $\Zm/(n+1)\Zm$. The action of $r$ on $\Aii(n)$ is determined by setting
$$r(c_n)=(-1)^nc_n$$
and using the rule
$$r\big((U,\sigma)\circ_i (V,\tau)\big)=\left\{\begin{array}{ll}
				r((U,\sigma))\circ_{i+1}(V,\tau),&i<k,\\
				&\\
				(-1)^{\mid U\mid \mid V\mid}r((V,\tau))\circ_1 r((U,\sigma)),&i=k,
				\end{array}\right.$$
where $U\in \Aii(k)$.

Geometrically, the operator $r$ acts on a tree $T$ by rotating it in the counterclockwise direction. That is, we relabel the leafs and root of $T$ so that the right most leaf of $T$ becomes the root of $r(T)$. See the following picture for an example:

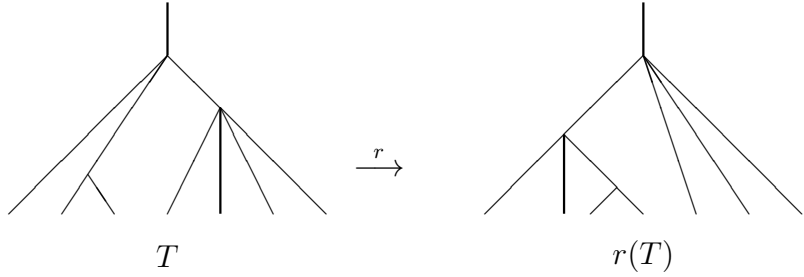
\begin{center}\label{T,rT}
 \setlength{\unitlength}{2pt}
\begin{picture}(150,50)(0,0)

\linethickness{0.3mm}
\put(30,40){\line(0,1){10}}
\linethickness{0.3mm}
\put(0,10){\line(1,1){30}}
\linethickness{0.3mm}
\put(10,10){\line(2,3){20}}
\linethickness{0.3mm}
\put(20,10){\line(-2,3){5}}
\linethickness{0.3mm}
\put(30,10){\line(1,2){10}}
\linethickness{0.3mm}
\put(40,10){\line(0,1){20}}
\linethickness{0.3mm}
\put(50,10){\line(-1,2){10}}
\linethickness{0.3mm}
\put(60,10){\line(-1,1){30}}

\linethickness{0.3mm}
\put(120,40){\line(0,1){10}}
\linethickness{0.3mm}
\put(90,10){\line(1,1){30}}
\linethickness{0.3mm}
\put(105,10){\line(0,1){15}}
\linethickness{0.3mm}
\put(110,10){\line(1,1){5}}
\linethickness{0.3mm}
\put(120,10){\line(-1,1){15}}
\linethickness{0.3mm}
\put(130,10){\line(-1,3){10}}
\linethickness{0.3mm}
\put(140,10){\line(-2,3){20}}
\linethickness{0.3mm}
\put(150,10){\line(-1,1){30}}

\put(30,2){\makebox(0,0)[cc]{$T$}}
\put(120,2){\makebox(0,0)[cc]{$r(T)$}}
\put(70,20){\makebox(0,0)[cc]{$\overset{r}{\longrightarrow}$}}

\end{picture}
\end{center}
By induction on the number of vertices, one can deduce the following formula.
\medskip
\begin{lem}~\label{rot-lemma}
Assume that $(T,\sigma)\in \Aii(n)$. Then we have
\[ r((T,\sigma))=(-1)^n (r(T),\sigma)\]
where $r(T)$, as in the previous example, is the ribbon tree obtained from $T$ by a counterclockwise rotation.
\end{lem}

\subsection{$A_\infty$-algebras}
Given a differential graded vector space $A$ we can define the operad $\End_A$ by taking
$$\End_A(n)=\Hom(A^{\otimes n}, A)$$
and defining the operadic composition and differential as follows:
\begin{align}
&f\circ_ig(a_1,\ldots,a_n)=(-1)^{\mid g\mid\sum_{l=1}^{i-1}\mid a_l\mid}f(a_1,\ldots,g(a_i,\ldots,a_{i+j-1}),\ldots,a_n),\nonumber\\
&\partial f(a_1,\ldots,a_n)=df(a_1,\ldots,a_n)-\sum_{i=1}^{n}(-1)^{\mid f\mid+\sum_{l=1}^{i-1}\mid a_l\mid}f(a_1,\ldots,da_i,\ldots,a_n),\nonumber
\end{align}
where $d$ is the differential on $A$.
\begin{defn}
An \emph{\Ai-algebra structure $(A, \rho)$} on $A$ is an operad homomorphism
$$\rho:\Aii\lto \End_A.$$
\end{defn}
Recall that $\Aii$ is generated as an operad by $c_j$, $j\geq 2$. Thus $\rho$ is equivalent to the choice of maps $\m_k: A^{\otimes k}\lto A$ of degree $k-2$, with $\m_1=d$, satisfying, for each $n$,
$$\sum_{i,j}(-1)^{i(j+1)+jn+j\sum_{l=1}^{i-1}\mid a_l\mid}\m_{n-j+1}(a_1,\ldots,\m_j(a_i,\ldots,a_{i+j-1}),\ldots,a_n)=0.$$
Here $\m_k=\rho(c_k)$, for $k\geq 2$ and the above equation is a consequence of the relation
$$\partial(c_n)=\sum_{i,j}(-1)^{i(j+1)+jn}c_{n-j+1}\circ_ic_j.$$

\subsection{Cyclic $A_\infty$-algebras}
Let $A$ be a finite dimensional differential graded vector space with a non-degenerate inner product $\langle ,\rangle$. Furthermore, we assume that the differential $d$ on $A$ is (graded) self-adjoint:
\[ \langle d a_0, a_1\rangle = (-1)^{\vert a_0 \vert + 1}\langle a_0, da_1\rangle.\]
This is equivalent to requiring that $\langle,\rangle: A\otimes A \lto \mQ$ is a chain map.
Then the endomorphism operad $\End_A$ is a cyclic operad (see \cite{GetKap}). Indeed, the inner product induces a natural identification
\[\Hom(A^{\otimes n}, A) \cong (A^\vee)^{\otimes n+1},\]
where $A^\vee = \Hom(A, \mQ)$ is the dual vector space of $A$.
Since the group $\ZZ/(n+1)\ZZ$ acts on the right hand side by permuting factors, we also get an action on $\Hom(A^{\otimes n}, A)$ via the above identification. 

A \emph{cyclic \Ai-algebra structure} on $A$ is then given by an operad homomorphism
\[\rho:\Aii\lto \End_A\] 
which is compatible with the cyclic structure. More explicitly, the following equation must hold for any $(T,\sigma)\in \Aii(n)$
$$\langle \rho(r((T,\sigma)))(a_0,\ldots,a_{n-1}),a_n\rangle =(-1)^{\mid a_0\mid\sum_{l=1}^n\mid a_l\mid}\langle \rho((T,\sigma))(a_1,\ldots,a_n),a_0\rangle .$$
Again, since $\Aii(n)$ is generated by corollas, this is equivalent to 
$$\langle \m_n(a_0,\ldots,a_{n-1}),a_n\rangle =(-1)^{n+\mid a_0\mid\sum_{l=1}^n\mid a_l\mid}\langle \m_n(a_1,\ldots,a_n),a_0\rangle .$$

\section{Tensor product of cyclic $A_\infty$-algebras}~\label{cyclic-sec}

In this section, we prove our first main result, Theorem~\ref{main-thm}. We will do this by showing the existence and uniqueness (up to homotopy) of a cyclic diagonal of the operad $\Aii$. 

\subsection{Diagonals of the $A_\infty$ operad}
Given two cyclic \Ai-algebras $(A,\rho_A,\langle ,\rangle_A)$ and $(B,\rho_B,\langle ,\rangle _B)$ we want to construct a cyclic \Ai-algebra on the differential vector space $(A\otimes B, d=d_A\otimes id+id\otimes d_B)$ equipped with the inner product
$$\langle a_1\otimes a_2,b_1\otimes b_2\rangle _{A \otimes B}:=(-1)^{\mid b_1\mid\mid a_2\mid}\langle a_1,b_1\rangle _A\langle a_2,b_2\rangle _B.$$

Let us recast this problem in operadic terms. First observe that given operads $\mO_1$ and $\mO_2$ we can define their tensor product as
$$(\mO_1\otimes\mO_2)(n)=\mO_1(n)\otimes\mO_2(n)$$
with the operadic composition and differential
\begin{align}
&(U_1\otimes U_2)\circ_i(V_1\otimes V_2)=(-1)^{\mid U_2\mid\mid V_1\mid}U_1\circ_i V_1\otimes U_2\circ_i V_2,\nonumber\\
&\partial(U_1\otimes U_2)=\partial U_1\otimes U_2+(-1)^{\mid U_1\mid} U_1\otimes \partial U_2.\nonumber
\end{align}
Also note that there is a map of operads 
$$i:\End_A\otimes \End_B\lto \End_{A\otimes B}$$
given by
$$i(f\otimes g)(a_1\otimes b_1,\ldots,a_n\otimes b_n)=(-1)^{\sum_{i < j}\mid b_i\mid\mid a_j\mid+\sum_{j=1}^n\mid a_j\mid \mid g\mid}f(a_1,\ldots,a_n)\otimes g(b_1,\ldots, b_n).$$

With these observations we see that if there exists a canonical (in an appropriate sense) map of operads
$$\Delta: \Aii\lto \Aii\otimes \Aii,$$
then we can simply use the composition
$$\rho=i\circ \rho_A\otimes\rho_B\circ\Delta: \Aii\lto \Aii\otimes\Aii \lto \End_A\otimes \End_B \lto \End_{A\otimes B}$$ 
to define the desired \Ai-algebra structure on $A\otimes B$.

Furthermore, if one can construct $\Delta$ which is also compatible with the cyclic group actions, then the resulting \Ai-algebra structure on $A\otimes B$ will also be cyclic. This follows from the observation that in order for the composition $i\circ \rho_A\otimes\rho_B\circ\Delta$ to be compatible with the cyclic structures, it is enough to require that $\Delta$ is so. This is because the morphisms $i$ and $\rho_A\otimes \rho_B$ are already compatible.

\subsection{ Existence of a cyclic diagonal} To characterize the morphism $\Delta$, we first recall some basic facts from~\cite{operads}. The operad $\Aii$ is a minimal resolution of the \emph{associative operad} $\As$, i.e. there is a morphism of dg operads (with $\As$ endowed with the trivial differential)
\[ \Aii\lto \As,\]
inducing an isomorphism on the homology operads $H_*(\Aii)\cong \As$. Furthermore, the operad $\As$ admits a canonical {\sl diagonal} morphism
\[ \As \lto \As\otimes \As\]
which is coassociative. In fact $\As$ is a {\sl Hopf operad} (see \cite{operads} for a definition).
\begin{defn}
A morphism of operads $\Delta: \Aii\lto \Aii\otimes \Aii$ is called a diagonal of $\Aii$ if the induced map on the corresponding homology operads
$$ H_*(\Aii)\lto H_*(\Aii\otimes\Aii)\cong H_*(\Aii)\otimes H_*(\Aii)$$
is the canonical diagonal of the associative operad $\As$ under the canonical isomorphism $H_*(\Aii)\cong \As$. If, in addition, the map $\Delta$ commutes with the cyclic action of $\ZZ/(n+1)\ZZ$, i.e.
$$\Delta\circ r= r\otimes r\circ \Delta,$$
then we say $\Delta$ is a \emph{cyclic diagonal} of $\Aii$.
\end{defn}

The existence of a diagonal (without the requirement of being cyclic) and its uniqueness up to homotopy follow from general results on Hopf operads. 
See \cite{operads} for a proof of this.
Here we give a different proof of the existence of a diagonal which will also be cyclic. The proof is based on the following fact:
\begin{prop}\label{contract}
The spaces $K(n)$ are contractible. Therefore
\begin{align}
H_*(\Aii(n))=H_*(C_*(K(n)))=\left\{\begin{array}{ll}
						\mQ&,*=0,\\
						0&,\textrm{otherwise}.
\end{array}
\right.\nonumber
\end{align}
\end{prop}
\proof For example, see~\cite{ShnSte} for a proof of the fact that each $K(n)$ can be realized as a convex polytope in $\R^N$, which implies the result.\qed

We will also need the following

\begin{lem}\label{alfa} Let $(C_*,\partial)$ be a chain complex with a $\Zm/(n+1)$ action (generated by $r$) and let $\alpha,\beta\in C_*$. The equation $\partial \alpha=\beta$ has a solution $\alpha$ satisfying $r(\alpha)=(-1)^n\alpha$ if and only if there exists $\gamma$ such that $\partial\gamma=\beta$ and $r(\beta)=(-1)^n\beta$.
\end{lem}
\proof
For the only if part, observe that 
$$r(\beta)=r(\partial\alpha)=\partial r(\alpha)=\partial((-1)^n\alpha)=(-1)^n\beta.$$
Conversely, we define
$$\alpha=\frac{1}{n+1}\sum_{i=0}^n(-1)^{in}r^i(\gamma)$$ and compute
\begin{align}
\partial\alpha&=\frac{1}{n+1}\sum_{i=0}^n(-1)^{in}\partial r^i(\gamma)=\frac{1}{n+1}\sum_{i=0}^n(-1)^{in}r^i(\partial\gamma)\nonumber\\
&=\frac{1}{n+1}\sum_{i=0}^n(-1)^{in}r^i(\beta)=\frac{1}{n+1}\sum_{i=0}^n(-1)^{in}(-1)^{in}\beta\nonumber\\
&=\beta.\nonumber
\end{align}
Moreover
\begin{align}
r(\alpha)&=\frac{1}{n+1}\sum_{i=0}^n(-1)^{in}r^{i+1}(\gamma)\nonumber\\
&= \frac{(-1)^n}{n+1}\sum_{i=0}^n(-1)^{(i+1)n}r^{i+1}(\gamma)\nonumber\\
&=\frac{(-1)^n}{n+1}\left(\sum_{j=1}^n(-1)^{jn}r^{j}(\gamma)+(-1)^{(n+1)n}r^{n+1}(\gamma)\right)\nonumber\\
&=\frac{(-1)^n}{n+1}\sum_{j=0}^n (-1)^{jn}r^{j}(\gamma)=(-1)^n
\alpha,\nonumber\end{align}
since $r^{n+1}=id$.\qed

\begin{thm}~\label{existence-thm}
There exists a cyclic diagonal for $\Aii$.
\end{thm}
\proof
The operad $\Aii$ is generated by the corollas $\{c_k\}_{k\geq 2}$. Therefore it is enough to define $\Delta(c_k)$ for $k\geq 2$. To ensure that $\Delta$ is a map of operads we just need to ensure that $\Delta$ is a chain map, that is, $\Delta(\partial c_k)=\partial\Delta(c_k)$. Also, to ensure that $\Delta$ respects the cyclic action, we only need to check it on corollas:
$$\Delta(r(c_k))=(r\otimes r)(\Delta(c_k)),\ k\geq 2.$$
We will proceed by induction on $k$. For $k=2$ we define
$$\Delta(c_2)=c_2\otimes c_2$$
and observe that it is a chain map since the differential is trivial on both sides. Moreover $r(c_2)=c_2$, which implies $\Delta$ commutes with $r$.

Now assume we have defined $\Delta(c_k)$ for all $k\leq n-1$ in such a way that it commutes with both $r$ and $\partial$. We want to define $\Delta(c_n)$ so that
$$\partial \Delta(c_n)=\Delta(\partial c_n)\textrm{ and }r\otimes r(\Delta(c_n))=\Delta(r(c_n)).$$
By Proposition \ref{contract}, $\Aii(n)$ is contractible, then so is $\Aii(n)\otimes \Aii(n)$ by the K\"unneth theorem. Therefore to prove that the first equation has a solution we only need to check that $\Delta(\partial c_n)$ is closed. By definition
\begin{align}
\Delta(\partial c_n)&=\Delta\left(\sum(-1)^{i(j+1)+jn}c_{n-j+1}\circ_ic_j\right)\nonumber\\
&=\sum(-1)^{i(j+1)+jn}\Delta(c_{n-j+1})\circ_i\Delta(c_j).\nonumber
\end{align}
By induction hypothesis the right hand side involves only $c_j$ with $j< n$. For these, $\Delta$ was already defined and commutes with $r$ and $\partial$. Therefore
$$\partial \Delta(\partial c_n)=\Delta(\partial ^2c_n)=0.$$
Hence we conclude that we can find $\gamma\in \Aii(n)\otimes\Aii(n)$ such that
$$\partial\gamma=\Delta(\partial c_n).$$
Now, in order to use Lemma \ref{alfa}, we need to compute
\begin{align}
r(\Delta(\partial c_n))&=r\left(\Delta\left(\sum(-1)^{i(j+1)+jn}c_{n-j+1}\circ_i  c_j\right)\right)\nonumber\\
&=\sum(-1)^{i(j+1)+jn}r\left(\Delta(c_{c-j+1})\circ_i \Delta(c_j)\right)\nonumber\\
&=\sum_j\Bigg(\sum_{1\leq i < n-j+1}(-1)^{i(j+1)+jn}r\left(\Delta(c_{n-j+1})\right)\circ_{i+1} \Delta(c_j)+\nonumber\\
&\hspace{2cm}+(-1)^{(n-j+1)(j+1)+jn+j(n-j+1)}r(\Delta(c_j))\circ_1 r(\Delta(c_{n-j+1}))\Bigg)\nonumber\\
&=\sum_{\substack{j\\2\leq i\leq n-j+1}}(-1)^{(i+1)(j+1)+jn+j+n+1}\Delta(c_{n-j+1})\circ_i\Delta(c_j)+\nonumber\\
&\hspace{2cm}+\sum_j(-1)^{(n+j+1)(j+1)+jn+j(n+j+1)+j+n+j+1}\Delta(c_j)\circ_1\Delta(c_{n-j+1})\nonumber
\end{align}
\begin{align}
&=(-1)^n\sum_{\substack{j\\2\leq i\leq n-j+1}}(-1)^{i(j+1)+jn}\Delta(c_{n-j+1})\circ_i\Delta(c_j)+ \nonumber\\
&\hspace{2cm}+\sum_k(-1)^{n+k+1+kn}\Delta(c_{n-k+1})\circ_1\Delta(c_k)\nonumber\\
&=(-1)^n\sum(-1)^{i(j+1)+jn}\Delta(c_{n-j+1})\circ_i\Delta(c_j)\nonumber\\
&=(-1)^{n}\Delta(\partial c_n).\nonumber
\end{align}
Lemma \ref{alfa} now implies that there is $\alpha\in\Aii(n)\otimes \Aii(n)$ satisfying 
$$\partial \alpha=\Delta(\partial c_n)\textrm{ and }r(\alpha)=(-1)^n\alpha.$$
Finally we define $\Delta(c_n)=\alpha$, which concludes the induction step.\qed

\begin{rmk} This proof can be immediately adapted to prove that we can choose a cocommutative diagonal. This means we can choose $\Delta$ so that 
$$\tau\circ\Delta=\Delta$$
where $\tau:\Aii(n)\otimes\Aii(n)\lto \Aii(n)\otimes\Aii(n)$ is the map that interchanges both factors of the tensor product.

Observe that Lemma \ref{alfa} works for the action of any finite group on a chain complex. Since the $\ZZ_2$-action generated by $\tau$ and the action generated by $r$ commute, we can in fact prove the existence of a diagonal which is simultaneously cyclic and cocommutative. With such diagonal, the tensor product of \Ai-algebras is naively commutative. In other words, given \Ai-algebras $A$ and $B$ the map $t=A\otimes B\lto B\otimes A$, $t(a\otimes b)=(-1)^{\vert a\vert\vert b\vert}b\otimes a$ is an isomorphism of \Ai-algebras. This is in contrast with the fact that we cannot construct a \emph{coassociative diagonal}, this meaning
$$(\Delta\otimes \text{id})\Delta\neq(\text{id}\otimes \Delta)\Delta,$$ 
for any choice of diagonal $\Delta$. See \cite{MarShn} for a proof of this fact.
\end{rmk}

\medskip
\subsection{Explicit formula for the diagonal}~\label{formula-sub}
As we mentioned in the Introduction there are three constructions of a diagonal on $\Aii$. The first is due to Saneblidze and Umble  \cite{SanUmb}, the second to Markl and Shnider \cite{MarShn} and the third to Loday \cite{Lod}. It seems these three constructions give the same diagonal, but this has not been checked for $\Delta(c_k)$, with $k$ larger than $5$.
 
 These diagonals are not cyclic. Here we will give explicit formulas for $\Delta(c_k)$, $k\leq 4$, which are cyclic. Unfortunately we do not know how to construct $\Delta(c_k)$ for $k\geq 6$ (we have constructed the case $k=5$ but will not present it in the present article due to its length). This appears to be a very complicated combinatorial problem.

\medskip
\noindent{\bf Diagonal for $k=2$.} We take $\Delta(c_2)=c_2\otimes c_2$.

\medskip
\noindent{\bf Diagonal for $k=3$.} If we denote by $B_1$ and $B_2$ the two trees with $3$ leaves and all internal vertices trivalent we have
$$\Delta(c_3)=\frac{1}{2}\left(B_1\otimes c_3+B_2\otimes c_3+c_3\otimes B_1+c_3\otimes B_2\right).$$ 
Note that this is also cocommutative.

\medskip
\noindent{\bf Diagonal for $k=4$.} Denote by $B_1$ trough $B_5$ the five trivalent trees with $4$ leaves with the orientation $e_1\wedge e_2$, where $e_1$ and $e_2$ are the internal edges indicated in the diagram below:
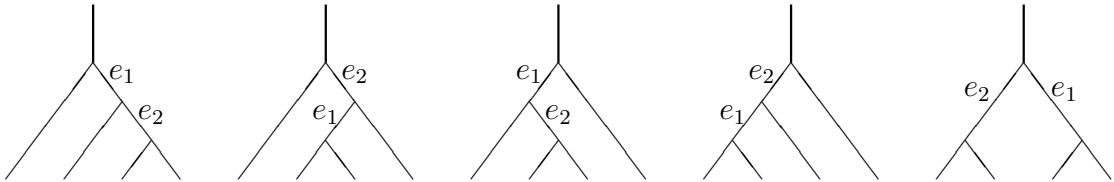
\begin{center}\label{k=4}
 \setlength{\unitlength}{2.2pt}
\begin{picture}(190,30)(0,0)

\linethickness{0.3mm}
\put(15,30){\line(0,-1){10}}
\linethickness{0.3mm}
\put(0,0){\line(3,4){15}}
\linethickness{0.3mm}
\put(30,0){\line(-3,4){15}}
\linethickness{0.3mm}
\put(10,0){\line(3,4){10}}
\linethickness{0.3mm}
\put(20,0){\line(3,4){5}}

\put(20,18){\makebox(0,0)[cc]{$e_1$}}
\put(25,11){\makebox(0,0)[cc]{$e_2$}}

\linethickness{0.3mm}
\put(55,30){\line(0,-1){10}}
\linethickness{0.3mm}
\put(40,0){\line(3,4){15}}
\linethickness{0.3mm}
\put(70,0){\line(-3,4){15}}
\linethickness{0.3mm}
\put(50,0){\line(3,4){10}}
\linethickness{0.3mm}
\put(60,0){\line(-3,4){5}}

\put(60,18){\makebox(0,0)[cc]{$e_2$}}
\put(55,11){\makebox(0,0)[cc]{$e_1$}}

\linethickness{0.3mm}
\put(95,30){\line(0,-1){10}}
\linethickness{0.3mm}
\put(80,0){\line(3,4){15}}
\linethickness{0.3mm}
\put(110,0){\line(-3,4){15}}
\linethickness{0.3mm}
\put(90,0){\line(3,4){5}}
\linethickness{0.3mm}
\put(100,0){\line(-3,4){10}}

\put(90,18){\makebox(0,0)[cc]{$e_1$}}
\put(95,11){\makebox(0,0)[cc]{$e_2$}}

\linethickness{0.3mm}
\put(135,30){\line(0,-1){10}}
\linethickness{0.3mm}
\put(120,0){\line(3,4){15}}
\linethickness{0.3mm}
\put(150,0){\line(-3,4){15}}
\linethickness{0.3mm}
\put(130,0){\line(-3,4){5}}
\linethickness{0.3mm}
\put(140,0){\line(-3,4){10}}

\put(130,18){\makebox(0,0)[cc]{$e_2$}}
\put(125,11){\makebox(0,0)[cc]{$e_1$}}

\linethickness{0.3mm}
\put(175,30){\line(0,-1){10}}
\linethickness{0.3mm}
\put(160,0){\line(3,4){15}}
\linethickness{0.3mm}
\put(190,0){\line(-3,4){15}}
\linethickness{0.3mm}
\put(170,0){\line(-3,4){5}}
\linethickness{0.3mm}
\put(180,0){\line(3,4){5}}

\put(182,15){\makebox(0,0)[cc]{$e_1$}}
\put(167,15){\makebox(0,0)[cc]{$e_2$}}

\end{picture}
\end{center}
Additionally we label the following trees as:
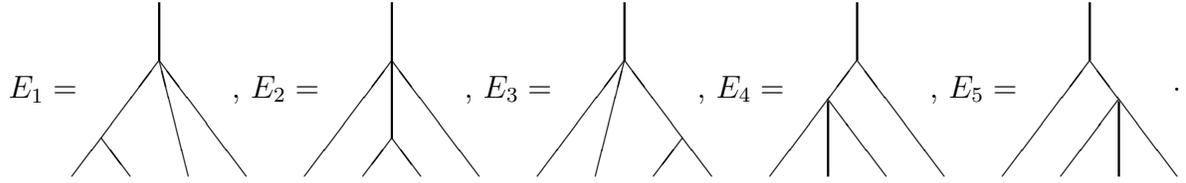
\begin{center}\label{Ei}
 \setlength{\unitlength}{2.2pt}
\begin{picture}(200,30)(-10,0)

\put(-5,15){\makebox(0,0)[cc]{$E_1=$}}

\linethickness{0.3mm}
\put(15,30){\line(0,-1){10}}
\linethickness{0.3mm}
\put(0,0){\line(3,4){15}}
\linethickness{0.3mm}
\put(30,0){\line(-3,4){15}}
\linethickness{0.3mm}
\put(10,0){\line(-3,4){5}}
\linethickness{0.3mm}
\put(20,0){\line(-1,4){5}}

\put(35,15){\makebox(0,0)[cc]{, $E_2=$}}

\linethickness{0.3mm}
\put(55,30){\line(0,-1){23.5}}
\linethickness{0.3mm}
\put(40,0){\line(3,4){15}}
\linethickness{0.3mm}
\put(70,0){\line(-3,4){15}}
\linethickness{0.3mm}
\put(50,0){\line(3,4){5}}
\linethickness{0.3mm}
\put(60,0){\line(-3,4){5}}

\put(75,15){\makebox(0,0)[cc]{, $E_3=$}}

\linethickness{0.3mm}
\put(95,30){\line(0,-1){10}}
\linethickness{0.3mm}
\put(80,0){\line(3,4){15}}
\linethickness{0.3mm}
\put(110,0){\line(-3,4){15}}
\linethickness{0.3mm}
\put(90,0){\line(1,4){5}}
\linethickness{0.3mm}
\put(100,0){\line(3,4){5}}

\put(115,15){\makebox(0,0)[cc]{, $E_4=$}}

\linethickness{0.3mm}
\put(135,30){\line(0,-1){10}}
\linethickness{0.3mm}
\put(120,0){\line(3,4){15}}
\linethickness{0.3mm}
\put(150,0){\line(-3,4){15}}
\linethickness{0.3mm}
\put(130,0){\line(0,1){13}}
\linethickness{0.3mm}
\put(140,0){\line(-3,4){10}}

\put(155,15){\makebox(0,0)[cc]{, $E_5=$}}

\linethickness{0.3mm}
\put(175,30){\line(0,-1){10}}
\linethickness{0.3mm}
\put(160,0){\line(3,4){15}}
\linethickness{0.3mm}
\put(190,0){\line(-3,4){15}}
\linethickness{0.3mm}
\put(170,0){\line(3,4){10}}
\linethickness{0.3mm}
\put(180,0){\line(0,1){13}}

\put(190,15){\makebox(0,0)[cc]{.}}

\end{picture}
\end{center}
Let $x\in \mQ$. The most general form of a cyclic $\Delta(c_4)$ is
$$\Delta(c_4)=\frac{1}{5}c_4\otimes \left(\sum_{i=1}^5 B_i\right)+\frac{1}{5}\left(\sum_{i=1}^5 B_i\right)\otimes c_4+\sum_{i,j}\gamma_{i,j} E_i\otimes E_j$$ with
\begin{align}
\left\{\begin{array}{l}
\gamma_{1,1}=-1/10-x,\\
\gamma_{1,2}=x,\\
\gamma_{1,3}=-2/5-x,\\
\gamma_{1,4}=-1/5+x,\\
\gamma_{1,5}=1/5+x,
\end{array}\right.\nonumber
\end{align}
together with the relations
\begin{align}
\left\{\begin{array}{l}
\gamma_{1,1}=\gamma_{2,2}=\gamma_{3,3}=\gamma_{4,4}=\gamma_{5,5},\\
\gamma_{1,2}=\gamma_{2,3}=\gamma_{3,4}=-\gamma_{4,5}=\gamma_{5,1},\\
\gamma_{1,3}=\gamma_{2,4}=-\gamma_{3,5}=-\gamma_{4,1}=\gamma_{5,2},\\
\gamma_{1,4}=-\gamma_{2,5}=-\gamma_{3,1}=-\gamma_{4,2}=\gamma_{5,3},\\
\gamma_{1,5}=\gamma_{2,1}=\gamma_{3,2}=\gamma_{4,3}=-\gamma_{5,4}.
\end{array}\right.\nonumber
\end{align}
Thus we see that the cyclicity condition uniquely determines $\Delta(c_3)$, but not $\Delta(c_4)$. However if we require that $\Delta$ must  also be cocommutative, it forces $x=-1/10$ and so it fixes $\Delta(c_4)$.

\medskip
\noindent{\bf Diagonal for $k=5$.} The most general form of $\Delta(c_5)$ has five degrees of freedom. Cocommutativity cuts down the ambiguity to four dimensions.

\subsection{Uniqueness up to cyclic homotopy} The goal of this subsection is to show that a cyclic diagonal of $\Aii$ is unique up to cyclic homotopy. On the level of algebras, this implies that the cyclic tensor product structure defined by any cyclic diagonal is unique up to a cyclic quasi-isomorphism.

We first recall the notion of a homotopy between morphisms of dg operads from~\cite{operads}. Denote by $\Omega^*_{[0,1]}$ the commutative dg algebra of algebraic differential forms on the unit interval $[0,1]$, with coefficients in $\mQ$, equipped with the standard de Rham differential. Note that since we are using the homological degree, we require the the space of one forms $\Omega^1_{[0,1]}$ to have degree $-1$. An elementary homotopy between two morphisms $f,g: \mO_1\lto \mO_2$ of dg operads is another homomorphism
\[ h: \mO_1\lto \mO_2\otimes_\mQ \Omega^*_{[0,1]}\]
such that $h(0)=f$ and $h(1)=g$~\footnote{A homotopy between $f$ and $g$ is a sequence of elementary homotopies $(h_1,\cdots, h_k)$ with $h_i(1)=h_{i+1}(0)$ where $1\leq i\leq k-1$, and $h_1(0)=f$ and $h_k(1)=g$. For purposes of this paper, the notion of elementary homotopy is enough.}.

For cyclic operads $\mO_1$, $\mO_2$, and cyclic morphisms $f$, $g$, we call an elementary homotopy $h$ {\sl cyclic} if it commutes with the cyclic action, where we put the trivial action on the $\Omega^*_{[0,1]}$ component.

\medskip
\begin{thm}~\label{uniquness-thm}
Let $\Delta_1, \Delta_2: \Aii\lto \Aii\otimes \Aii$ be two cyclic diagonals, then there exists a cyclic homotopy $h: \Aii \lto \Aii\otimes\Aii\otimes \Omega^*_{[0,1]}$ between them.
\end{thm}

\proof The proof is very similar to that of Theorem~\ref{existence-thm}. We shall only sketch it here. Indeed, the idea is to use the acyclicity of $ \Aii\otimes\Aii\otimes \Omega^*_{[0,1]}$ (which follows from acyclicity of $\Aii$ and $\Omega^*_{[0,1]}$) to inductively construct $h$. The only difference is that when constructing $h(c_n)$, we also need to impose the boundary condition
\[ h(c_n)(0)=\Delta_1(c_n), \mbox{ and } h(c_n)(1)=\Delta_2(c_n).\]
For this we argue as follows. For $n=2$ we simply take $h(c_2)=\Delta_1(c_2)=\Delta_2(c_2)$. As in the proof of Theorem~\ref{existence-thm} we can find $h'(c_n)$, some extension of $h(c_2),\cdots, h(c_{n-1})$ as a morphism of dg operads (ignoring cyclicity and the boundary condition). Since $h(c_j)$ $(2\leq j\leq n-1)$ has the required boundary condition, we conclude that
\begin{align*}
 \partial \big( h'(c_n)(0) - \Delta_1(c_n)\big)&=0,\\
  \partial \big(h'(c_n)(1) - \Delta_2(c_n)\big)&=0.
\end{align*}
Thus by acyclicity, there exist elements $a, b\in (\Aii\otimes\Aii)(n)$ of degree $n-1$ such that
\begin{align*}
\partial a &=h'(c_n)(0)-\Delta_1(c_n), \\
\partial b &=h'(c_n)(1)-\Delta_2(c_n). 
\end{align*}
Next we set 
$$h''(c_n)=h'(c_n)-[(1-t)\partial a+t \partial b]-(-1)^{n}(a-b)dt,$$ 
and, finally, take $h(c_n)$ to be the average of $h''(c_n)$ by the cyclic group $\ZZ/(n+1)\ZZ$. Note that $h(c_n)$ gives a cyclic dg-operad map with the correct boundary condition. \qed 

On the level of algebras, we can see that the cyclic $A_\infty$-algebra structure on $A\otimes B$ defined using a cyclic diagonal of $\Aii$ is unique up to cyclic quasi-isomorphisms. Indeed, a cyclic homotopy $h$ induces a morphism of dg operads defined by the composition
\[   (i\otimes\id) (\rho_1\otimes\rho_2\otimes\id)\circ h: \Aii \rightarrow \Aii\otimes\Aii\otimes \Omega^*_{[0,1]} \rightarrow \End_{A}\otimes \End_B \otimes \Omega^*_{[0,1]} \rightarrow \End_{A\otimes B}\otimes \Omega^*_{[0,1]} ,\]
which we denote by $\theta$. Using the isomorphism $\Aii\cong \Omega\As^{\mbox{\textexclamdown}}$ where $\As^{\mbox{\textexclamdown}}$ is the Koszul dual cooperad of $\As$ and $\Omega$ is the coBar construction, we conclude that
\begin{align*}
\Hom_{\mbox{dg. op.}}(\Aii, \End_{A\otimes B}\otimes \Omega^*_{[0,1]}) &\cong \Tw( \As^{\mbox{\textexclamdown}}, \End_{A\otimes B}\otimes\Omega^*_{[0,1]})\\
&\cong \Hom_{\mbox{dg. coaug. coop.}}(\As^{\mbox{\textexclamdown}}, \BB \End_{A\otimes B}\otimes\Omega^*_{[0,1]})\\
&\cong \Codiff\big(\As^{\mbox{\textexclamdown}}(\End_{A\otimes B})\otimes \Omega^*_{[0,1]}\big).
\end{align*}
Here the functor $\BB$ is the Bar construction, which is the right adjoint of the coBar functor $\Omega$. Via the above isomorphism, the morphism $\theta$ corresponds to an element of the set $\Codiff\big(\As^{\mbox{\textexclamdown}}(\End_{A\otimes B})\otimes \Omega^*_{[0,1]}\big)$ which is compatible with the inner product on $A\otimes B$. In~\cite{Fukaya}, such an element was called a {\sl cyclic pseudo-isotopy} between the two cyclic \Ai-structures on $A\otimes B$ defined using $\Delta_1$ and $\Delta_2$. It was proved in Proposition $9.2$ of \cite{Fukaya} that cyclic pseudo-isotopies give rise to cyclic quasi-isomorphisms between the boundary cyclic $A_\infty$-structures defined by $\theta(0)$ and $\theta(1)$.

\section{Tensor products and cup products}~\label{graph-sec}

In this section, we prove our second main result, Theorem~\ref{main2-thm} in the case of two even \Ai-algebras. We leave the remaining cases for the next section. The bulk of the work is to use a cyclic diagonal of $\Aii$ to construct a diagonal of the graph complex $\mG_*$.

\subsection{The graph homology complex}

We recall the definition of the graph homology complex $\mG_*$, following Igusa's careful treatment of signs and orientations~\cite{Igu}.

\begin{defn}
A \emph{ribbon graph} is a (finite) connected graph for which every vertex have valency greater or equal than $3$, together with a fixed cyclic ordering on the half-edges incident to each vertex.

An \emph{orientation of a ribbon graph} $\Gamma$ is an integral generator of $\det\big(V(\Gamma)\oplus H(\Gamma)\big)$, where $V(\Gamma)$ is the set of vertices, $H(\Gamma)$ the set of half edges of $\Gamma$ and $\det$ is the top exterior power. In other words, it is an orientation on the vector space spanned by the vertices and half-edges of $\Gamma$.
\end{defn}
We define $\mG_*$ to be the vector space generated by isomorphism classes of pairs $( \Gamma,\sigma) $ where $\Gamma$ is a ribbon graph and $\sigma$ is an orientation of $\Gamma$, with the relation $(\Gamma,-\sigma) =-(\Gamma,\sigma)$. Here an isomorphism between $(\Gamma, \sigma) $ and $( \Gamma',\sigma') $ is an isomorphism $\varphi:\Gamma\lto\Gamma'$  of ribbon graphs which also preserves orientation.

The grading on $\mG_*$ is defined as follows. Given a pair $(\Gamma,\sigma)$, we denote by $v_1,\ldots,v_k$ the vertices of $\Gamma$, and define 
$$\vert(\Gamma,\sigma)\vert=\sum_{i=1}^{k}(\val(v_i)-3)$$
where $\val(v_i)$ is the valency of the vertex $v_i$.

To define the boundary map on $\mG_*$, we need to introduce the {\sl contraction} and {\sl expansion} operations. Given an edge $e$ of a ribbon graph $\Gamma'$ which is not a \emph{loop}, we define another ribbon graph $\Gamma:=\Gamma'/e$ to be the ribbon graph obtained from $\Gamma'$ by contracting the edge $e$. Let us denote by $v^-$ and $v^+$ the boundary vertices of $e$, and $e^-$, $e^+$ the two half edges of $e$ incident to $v^-$ and $v^+$ respectively. Denote by $v$ the vertex of $\Gamma$ obtained from identifying $v^-$ and $v^+$. Assume that the cyclic orderings of the half-edges at $v^-$ and $v^+$ were:
\[ e_1\cdots e_ke^-, \mbox{ and } e^+f_1\cdots f_n.\]
We define the cyclic ordering at the vertex $v$ of $\Gamma$ by 
$$e_1\ldots e_kf_1\ldots f_n.$$ 
Also, given an orientation $\sigma'$ of $\Gamma'$ of the form $\epsilon\langle v^-v^+e^-e^+\ldots\rangle $ where $\epsilon=\pm 1$~\footnote{This can always be achieved by a permutation action on $\sigma'$, then put $\epsilon$ to be the sign of this permutation.}, we define a natural orientation on $\Gamma=\Gamma'/e$ by setting $\sigma'/e:=\epsilon\langle v\ldots\rangle $. 

In the reverse direction, given a vertex $v$ of a ribbon graph $\Gamma$, one can expand $v$ into two vertices $v^-$, $v^+$ together with an edge $e$ joining them, in $\frac{\val(v)^2-3\val(v)}{2}$ different ways. This way we obtain a new oriented ribbon graph $(\Gamma',\sigma')$ with a distinguished edge $e$, satisfying $(\Gamma'/e,\sigma'/e)=(\Gamma, \sigma)$. 

Finally we define
$$\partial(\Gamma,\sigma)=\sum (\Gamma',\sigma')$$
where the sum is taken over pairs consisting a vertex $v$ of $\Gamma$ and a choice of expansion at $v$. 

\begin{defn}We define the ribbon graph complex to be the pair $(\mG_*,\partial)$.
\end{defn}

\subsection{Diagonals of the graph complex} In this subsection we will first construct a diagonal in $\mG_*$, that is a chain map 
\[ \delta: \mG_*\rightarrow \mG_*\otimes \mG_*.\]
We will then show that we can use this map to compute the cup product in the cohomology of the moduli space of Riemann surfaces by comparing $\delta$ to the Alexander-Whitney diagonal on a simplicial model for this space, namely the simplicial complex of the nerve of the category $\mathcal{F}at$ defined by Igusa in \cite{Igu}. This category has as objects ribbon graphs and as morphisms contraction of edges between ribbon graphs. Igusa proves in~\cite[Theorem 1.22]{Igu}, that there is a rational homotopy equivalence
\[ \phi: C_*(\mathcal{F}at) \rightarrow \mG_*,\]
where $C_*(\mathcal{F}at)$ is the chain complex of the nerve of $\mathcal{F}at$. Since the nerve of $\mathcal{F}at$ is a simplicial complex, the complex $C_*(\mathcal{F}at)$ admits a canonical diagonal: the Alexander-Whitney diagonal. In Proposition \ref{diagonal-prop} we will show that $\delta$ and the Alexander-Whitney diagonal of $C_*(\mathcal{F}at)$ are homotopic via the equivalence $\phi$.

We start by recalling some facts about $\phi$. The map $\phi$ is constructed as the (unique up to homotopy) chain map carried by an acyclic carrier called the ``forest carrier" $F_*$, introduced by Igusa in \cite[Section 1.5]{Igu}. Recall such an acyclic carrier is simply a functor from $\mathcal{F}at$ to the category of augmented, acyclic chain complexes over $\mG_*$. Explicitly, the carrier $F_*$ consists of the following data 
\begin{itemize}
\item[(a)] an augmented and acyclic chain complex $F_*(\Gamma)$ for each object $\Gamma\in \mathcal{F}at$;
\item[(b)] an augmented chain map
\[ p_\Gamma: F_*(\Gamma) \rightarrow \mG_*;\]
\item[(c)] a chain map 
\[ p_{(\Gamma'\rightarrow\Gamma)}: F_*(\Gamma')\rightarrow F_*(\Gamma),\]
for each morphism $(\Gamma'\rightarrow\Gamma)$.
\end{itemize}
These maps are required to be natural, which means they satisfy
\begin{align*}
p_\Gamma\circ p_{(\Gamma'\rightarrow\Gamma)}&=p_{\Gamma'}\\
p_{(\Gamma'\rightarrow\Gamma)}\circ p_{(\Gamma''\rightarrow\Gamma')}&=p_{(\Gamma''\rightarrow\Gamma)}.
\end{align*}

The forest carrier is then defined as follows: $F_*(\Gamma)$ is the vector space spanned by isomorphism classes of oriented ribbon graphs over $\Gamma$, that is morphisms $\big(\Gamma_1\rightarrow \Gamma\big)$ together with an orientation $\sigma_1$ on $\Gamma_1$. Two such objects are isomorphic if there is an orientation preserving isomorphism $\Gamma_1 \rightarrow \Gamma_2$ which intertwines the maps to $\Gamma$. To define the augmentation, recall from \cite{Igu} that a ribbon graph $\Gamma_0$ of degree zero has a canonical orientation $c(\Gamma_0)$. Then we define the augmentation as follows: $\epsilon(\Gamma_0, \sigma_0)=1$ if $\sigma_0$ agrees with $c(\Gamma_0)$ and $-1$ otherwise. The morphism of augmented chain complexes
\[ p_\Gamma: F_*(\Gamma) \rightarrow \mG_*\]
is defined by sending $\big(\Gamma'\rightarrow \Gamma\big)$ to $\Gamma'$ and the maps $p_{(\Gamma_1\rightarrow\Gamma_2)}$ are given by by post-composing with the map $\Gamma_1\rightarrow \Gamma_2$. We shall often omit the orientation $\sigma'$ from the notations.

As explained in \cite[Section 2]{Igu2} an acyclic carrier, carries a unique chain map 
\[ \phi: C_*(\mathcal{F}at) \rightarrow \mG_*.\]
This means that for any generator $(\Gamma_0\rightarrow \ldots \rightarrow \Gamma_n) \in C_*(\mathcal{F}at)$, there is $\tilde{\phi}(\Gamma_0\rightarrow \ldots \rightarrow \Gamma_n) \in F_*(\Gamma_n)$ such that $\phi(\Gamma_0\rightarrow \ldots \rightarrow \Gamma_n)=p_{\Gamma_n}(\tilde{\phi}(\Gamma_0\rightarrow \ldots \rightarrow \Gamma_n) \in F_*(\Gamma_n)$. Moreover this is supposed to respect the augmentations and the boundary maps, see \cite{Igu, Igu2} for full details.

\vspace{.1cm}

We are now ready to define $\delta$. This will depend on the choice of a cyclic diagonal of $\Aii$ constructed in the previous section, in the following we shall fix, once and for all, such a diagonal $\Delta$. 
We first introduce the following

\medskip
\begin{constr}~\label{delta-constr}
Let $(\Gamma,\sigma)\in\mG_*$, pick an ordering of the vertices of $\Gamma$, $v_1,\ldots,v_k$ and pick a half-edge incident at each vertex $v_i$ and label it by $e_{i,0}$. Then label the other half-edges following the cyclic ordering $e_{i,0},e_{i,1},\ldots, e_{i,n_i}$ (if $\val(v_i)=n_i+1$). We define $A_\Gamma\in\{\pm1\}$ by the following equality
$$\sigma=A_\Gamma \ v_1\wedge e_{1,0}\wedge e_{1,1}\ldots e_{1,n_1}\wedge v_2 \wedge e_{2,0} \wedge e_{2,1}\ldots e_{2,n_2}\ldots v_k \wedge e_{k,0}\wedge e_{k,1}\ldots e_{k,n_k} .$$

Now, given oriented trees $(T_i,\sigma_i=f_{i,1}\wedge\ldots\wedge f_{i,t_i})\in C_{n_i-t_i-2}(K(n_i))$ with $n_i$ leaves for $i=1,\ldots,k$ we define the ribbon graph $\Gamma(T_1,\ldots,T_k)$ by replacing each vertex $v_i$ by the tree $T_i$, gluing the root of $T_i$ to $e_{i,0}$ and the leaves of $T_i$ to $e_{i,1}, \ldots, e_{i,n_i}$, following the cyclic ordering. We define the orientation $\sigma_{\Gamma(T_1,\ldots,T_k)}$ on $\Gamma(T_1,\ldots,T_k)$ so that
$$\sigma_{\Gamma(T_1\ldots T_k)}/f_{1,1}/\ldots/f_{1,t_1}/f_{2,t_1}\ldots f_{2,t_2}/\ldots/f_{k,1}\ldots/f_{k,t_k}=\sigma.$$
Note that there is a natural map of ribbon graphs $\big(\Gamma(T_1\ldots T_k\big)\rightarrow \Gamma)$ which contracts internal edges in $T_i$ to the vertex $v_i$ of $\Gamma$.
\end{constr}

\begin{defn}\label{def-delta}
For a cyclic diagonal $\Delta$ of $\Aii$, we use the Sweedler notation
$$\Delta(c_n):= c_n^{(1)} \otimes c_n^{(2)}$$
to denote the image of $\Delta$ applied to the $n$-th corolla, and we denote by $t_n^{(i)}$ $(i=1,2)$ the number of internal edges of $c_n^{(i)}$. Define a linear map $\delta: \mG_*\rightarrow\mG_*\otimes\mG_*$ by formula:
\[ \delta(\Gamma)=(-1)^{|\Gamma|}A_\Gamma (-1)^{\sum_{i<j} n_it_{n_j}^{(1)}+t_{n_i}^{(2)}t_{n_j}^{(2)}}\Gamma(c_{n_1}^{(1)},\ldots,c_{n_k}^{(1)})\otimes \Gamma (c_{n_1}^{(2)},\ldots,c_{n_k}^{(2)}).\]
\end{defn}

We will now prove several preliminary lemmas about $\delta$ and Construction~\ref{delta-constr}.

\begin{lem}
The definition of $\delta$ is independent of the choices made in Construction~\ref{delta-constr}.
\end{lem}

\proof If we interchange the two vertices $v_1$ and $v_2$, the sign $A_\Gamma$ changes by $(-1)^{n_1n_2}$. The sign $(-1)^{\sum_{i<j} n_it_{n_j}^{(1)}+t_{n_i}^{(2)}t_{n_j}^{(2)}}$ changes by
$(-1)^{n_1t_{n_2}^{(1)}+n_2t_{n_1}^{(1)}}$. Furthermore, the orientation of $\Gamma(c_{n_1}^{(1)},\ldots,c_{n_k}^{(1)})\otimes \Gamma (c_{n_1}^{(2)},\ldots,c_{n_k}^{(2)})$ changes by a factor $(-1)^{t_{n_1}^{(1)}t_{n_2}^{(1)}+t_{n_1}^{(2)}t_{n_2}^{(2)}}$. Thus the total sign change is given by $(-1)$ to the power of
\begin{align*}
n_1n_2 &+n_1t_{n_2}^{(1)}+n_2t_{n_1}^{(1)}+t_{n_1}^{(1)}t_{n_2}^{(1)}+t_{n_1}^{(2)}t_{n_2}^{(2)}\\
&= (n_1+t_{n_1}^{(1)})(n_2+t_{n_2}^{(1)})+t_{n_1}^{(2)}t_{n_2}^{(2)}\\
&= 0.
\end{align*}
Note that the last equality follows from the fact that $t_{n_i}^{(1)}+t_{n_i}^{(2)}=n_i \pmod 2$.

Changing the half-edge $e_{i,0}$ to $e_{i,1}$ corresponds to rotating $\Delta(c_n)$. By cyclicity of the diagonal and Lemma \ref{rot-lemma} this changes the sign by $(-1)^{n_i}$ which cancels with the change in $A_\Gamma$.\qed

\begin{lem}\label{lema1}
Let $(\Gamma,\sigma)$ and $T_i$ ($i=1,\ldots,k$) be as in Construction~\ref{delta-constr}. Then we have
$$\partial\left(\Gamma(T_1,\ldots,T_k)\right)=\sum_{i=1}^k(-1)^{\sum_{j< m}t_j}\Gamma(T_1,\ldots,\partial T_m,\ldots, T_k).$$
Here $t_j$ is the number of internal edges of $T_j$.
\end{lem}
\proof
Assume that
$$\partial T_m=\sum(U,\sigma_k=e\wedge \sigma_{T_m})$$ 
where $U/e\simeq T_m$. Then the orientation $\sigma_{T_1,\ldots,U,\ldots,T_k}$ on $\Gamma(T_1,\ldots,T_{m-1},U,T_{m+1},\ldots,T_k)$ is determined by
$$\sigma_{T_1,\ldots,U,\ldots,T_k}/\sigma_1/\ldots/e/\sigma_{T_m}/\ldots/\sigma_k=\sigma.$$
Then we have
$$(-1)^{\sum_{j< m}t_j}\sigma_{T_1,\ldots,U,\ldots,T_k}/e/\sigma_1/\ldots/\sigma_k=\sigma,$$
which implies the lemma.\qed

The following is a reformulation of a result of Igusa proving the acyclicity of the carrier $F_*$.

\medskip
\begin{cor}~\label{ch-map-cor}
Let $F_*$ denote Igusa's forest carrier. Furthermore, choose an orientation $\sigma$ on $\Gamma\in \mathcal{F}at$. There is an isomorphism of chain complexes
$$\phi_\Gamma: C_*(K(n_1))\otimes \ldots \otimes C_*(K(n_k)) \rightarrow F_*(\Gamma)$$
defined by 
$$\phi_\Gamma(T_1,\ldots,T_k):=(-1)^{|\Gamma|}A_\Gamma (-1)^{\sum_{i<j} n_it_j} \big(\Gamma(T_1,\ldots,T_k)\rightarrow \Gamma\big).$$
The definition of $\phi_\Gamma$ is independent of the choice of $\sigma$.
\end{cor} 
\proof  The sign follows from Lemma~\ref{lema1}. That $\phi_\Gamma$ is an isomorphism is proved by Igusa~\cite[Proposition 1.21]{Igu}. If we change the orientation of $\Gamma$, the orientation of $\Gamma(T_1,\ldots,T_k)$ also changes and $A_\Gamma$ gets multiplied by $-1$. As in the definition of $F_*(\Gamma)$ we equated reversal of orientation with change of sign, these cancel. \qed

\begin{defn}~\label{delta-gamma-def}
Define a morphism $\delta_\Gamma: F_*(\Gamma) \rightarrow F_*(\Gamma)\otimes F_*(\Gamma)$ as the composition of the maps in the following diagram.
\[\begin{CD}
F_*(\Gamma) @>\phi_\Gamma^{-1} >>     C_*(K(n_1))\otimes\ldots\otimes C_*(K(n_k)) \\
        @.                                      @VV\Delta\otimes\ldots\otimes\Delta V\\
         @. \big(C_*(K(n_1))\otimes C_*(K(n_1))\big)\otimes\ldots\otimes \big(C_*(K(n_k))\otimes C_*(K(n_k))\big)\\
          @.                                       @VV\tau V\\
         F_*(\Gamma)\otimes F_*(\Gamma) @<\phi_\Gamma\otimes\phi_\Gamma<<   \big(C_*(K(n_1))\otimes\ldots\otimes  C_*(K(n_k))\big)\otimes\big(C_*(K(n_1))\otimes\ldots\otimes C_*(K(n_k))\big)\\
\end{CD}\]
Here $\tau$ is the permutation of tensor components as indicated in the diagram. It follows from Corollary~\ref{ch-map-cor} that $\delta_\Gamma$ is a chain map.
\end{defn}

\begin{lem}\label{lema.delta}
The map $\delta_\Gamma$ does not depend on the choice of the ordering of $v_i$'s, nor on the choices of roots $e_{i,0}$'s.
\end{lem}

\proof First we introduce some notation. As in Definition \ref{def-delta}, let $\Delta(T_i)=T_i^{(1)} \otimes T_i^{(2)}$ and denote by $t_i^{(l)}$ the number of internal edges of $T_i^{(l)}$. Now suppose we interchange two consecutive vertices $v_1$ and $v_2$. Then $A_\Gamma$ changes by a factor of $(-1)^{n_1n_2}$. One easy computation shows that the contributions of the factor $(-1)^{\sum_{i<j} n_it_j}$ in the three maps $\phi_\Gamma$ (or its inverse) in the definition of $\delta_\Gamma$ cancel out. From Construction~\ref{delta-constr}, we can see that the orientation on $\Gamma(T_1,\ldots,T_k)$ changes by $(-1)^{t_1 t_2}$ and the orientation on $\Gamma(T_1^{l},\ldots,T_k^{l})$ changes by $(-1)^{t_1^{l} t_2^{l}}$, for $l=1,2$. Finally the map $\tau$ contributes with an extra sign of 
$$(-1)^{(n_1+t_1^{(2)})(n_2+t_2^{(1)})+(n_1+t_1^{(1)})(n_2+t_2^{(2)})}=(-1)^{n_1 t_2+n_2 t_1+t_1^{(2)}t_2^{(1)}+t_1^{(1)}t_2^{(2)}}, $$ 
where we have used that $n_i + t_i=t_i^{(1)}+t_i^{(2)} \pmod 2$, which follows from the fact that $\Delta$ is degree $0$.
Putting these signs together
\begin{align*}
n_1n_2 +  t_1t_2 &+t_1^{(1)}t_2^{(1)}+ t_1^{(2)}t_2^{(2)} + n_1 t_2+n_2 t_1+t_1^{(2)}t_2^{(1)}+t_1^{(1)}t_2^{(2)}=\\
&=n_1n_2 + n_1 t_2+n_2 t_1 + t_1t_2 +(t_1^{(1)}+t_1^{(2)})(t_2^{(1)}+t_2^{(2)})=\\
&=n_1n_2 + n_1 t_2+n_2 t_1 + t_1t_2 +(n_1 +t_1)(n_2+t_2)=0 \pmod 2,
\end{align*}
we conclude that $\delta_\Gamma$ does not depend on the order of vertices.

Next suppose we change the root $e_{1,0}$ by cyclically permuting the half-edges around $v_1$,
$$e_{1,0}e_{1,1}\ldots e_{1,n_1}\lto e_{1,n_1}e_{1,0}e_{1,1}\ldots e_{1,n_1-1}.$$
Then $A_\Gamma$ changes by $(-1)^{n_1}$. By Lemma~\ref{rot-lemma} this extra sign also appears from the rotation action on $T_1$. This together with the cyclicity of $\Delta$ implies the result.\qed

Let $\Gamma'\rightarrow \Gamma$ be a morphism in $\mathcal{F}at$ which contracts only one non-loop edge $e$. Furthermore, we choose orientations $\sigma'$ and $\sigma$ requiring that $(\Gamma'/e,\sigma'/e)=(\Gamma,\sigma)$. We introduce some notations. Assume $\Gamma'$ comes from expanding $\Gamma$ at the vertex $v_m$. We order vertices of $\Gamma'$ by 
$$v_1\ldots v_{m-1}v^+_mv^-_mv_{m+1}\ldots v_k.$$ 
where $v_m^+$ in the unique vertex that contains the half-edge $e_{m,0}$ after expansion. We then order half-edges of $v^+$ and $v^-$ by 
\begin{align*}
e_{m,0}&e_{m,1}\ldots e_{m,i-1}e^+e_{m,i+j}\ldots e_{m,n_m} \mbox{ and }\\
e^-&e_{m,i}\ldots e_{m,i+j-1}
\end{align*}
respectively. Here $e^-$ and $e^+$ are the new half-edges incident to $v_m^-$ and $v_m^+$. We use this choice of ordering to write down the maps $\delta_\Gamma$ and $\delta_{\Gamma'}$.

\begin{prop}\label{compatible-prop}
Let $\Gamma$ and $\Gamma'$ be as above. Then the following diagram is commutative.
\[\begin{CD}
F_*(\Gamma') @>>> F_*(\Gamma)\\
@V \delta_{\Gamma'} VV             @V\delta_{\Gamma} VV\\
F_*(\Gamma')\otimes F_*(\Gamma') @>>> F_*(\Gamma)\otimes F_*(\Gamma)
\end{CD}\]
Here horizontal arrows are the maps induced by the morphism $\Gamma' \rightarrow\Gamma$.
\end{prop}

The proposition is clear on the level of ribbon graphs. So the content of it is to match the orientations. For this we first need to prove two lemmas.

\begin{lem}\label{lema2}
Given $U\in C_*(K(n_m-j+1))$ and $V\in  C_{j-t_V-2}(K(j))$ (with  $t_V$ the number of internal edges of $V$) we have
$$\Gamma'(T_1,\ldots,T_{m-1},U,V,T_{m+1},\ldots, T_k)=(-1)^\epsilon \ \Gamma(T_1,\ldots,T_{m-1},U\circ_iV,T_{m+1},\ldots, T_k)$$
with $\epsilon=i(j+1)+(j+t_V)(n_m-j+1)+\sum_{l> m}t_l$. Additionally,
$$A_{\Gamma'}=(-1)^{\sum_{l< m}n_l+i(j+1)+n_mj+1}A_\Gamma.$$
\end{lem}
\proof
The equality as ribbon graphs is clear, we just need to compare orientations. Denote by $\sigma_U$ and $\sigma_V$ the orientations of $U$ and $V$. By definition of composition we have
$$\sigma_{U\circ_iV}=(-1)^{i(j+1)+(j+t_V)(n_m-j+1)}\sigma_U\wedge\sigma_V\wedge e.$$
Also by definition,
$$\sigma'_{T_1,\ldots,U,V,\ldots, T_k}/\sigma_1/\ldots/\sigma_U/\sigma_V/\ldots\sigma_k=\sigma',$$
$$\sigma'/e=\sigma,$$
and
$$\sigma_{T_1,\ldots,U\circ_iV,\ldots, T_k}/\sigma_1/\ldots/\sigma_{U\circ_i V}/\ldots/\sigma_k=\sigma.$$
Then
\begin{align}&(-1)^{i(j+1)+(j+t_V)(n_m-j+1)}\sigma_{T_1,\ldots,U\circ_iV,\ldots, T_k}/\sigma_1/\ldots/\sigma_U/\sigma_V/e/\ldots\sigma_k=\sigma\nonumber\\
\Rightarrow \ \ & (-1)^{i(j+1)+(j+t_V)(n_m-j+1)+\sum_{l> m}t_l}\sigma_{T_1,\ldots,U\circ_iV,\ldots, T_k}/\sigma_1/\ldots/\sigma_U/\sigma_V/\ldots/\sigma_k/e=\sigma,\nonumber
\end{align}
which implies
$$\sigma'_{T_1,\ldots,U,V,\ldots, T_k}=(-1)^{i(j+1)+(j+t_V)(n_m-j+1)+\sum_{l> m}t_l}\sigma_{T_1,\ldots,U\circ_iV,\ldots, T_k}.$$
For the second statement note that, by definition of $(\Gamma',\sigma')$ we have $\sigma'=v^+_m\wedge e^+\wedge e^-\wedge\sigma$, if we identify the vertex $v_m^-$ in $\Gamma'$ with $v_m$ in $\Gamma$. Also by definition of $A$, we have:
\begin{align}
\sigma'&=A_{\Gamma'}\langle v_1e_{1,0}\ldots e_{1,k}\ldots v_m^+e_{m,0}\ldots e_{m,i-1}e^+e_{m,i+j}\ldots e_{m,n_m}v^-_me^-e_{m,i}\ldots e_{m,i+j-1}\ldots\rangle \nonumber\\
&=(-1)^{i(j+1)+jn_m+1}A_{\Gamma'}\langle v_1\ldots v^+_mv_m^-e^+e^-e_{m,0}\ldots e_{m,n_m}\ldots\rangle \nonumber\\
&=(-1)^{i(j+1)+jn_m+1}A_{\Gamma'}(-1)^{\sum_{l< m}n_l}\langle v_m^+e^+e^-v_1\ldots v_m^-e_{m,0}\ldots e_{m,n_m}\ldots\rangle \nonumber\\
&=(-1)^{i(j+1)+jn_m+1+\sum_{l< m}n_l}A_{\Gamma'}\ v_m^+\wedge e^+\wedge e^-\wedge(A_{\Gamma}\sigma).\nonumber
\end{align}
Hence
$$A_{\Gamma'}=A_\Gamma(-1)^{i(j+1)+jn_m+1+\sum_{l< m}n_l},$$
which completes the proof of the lemma.\qed

\begin{lem}\label{inclusion-lem}
Let $\Gamma$ and $\Gamma'$ be as in Proposition~\ref{compatible-prop}. Then the following diagram is commutative.
\[\begin{CD} 
F_*(\Gamma')@>\phi_{\Gamma'}^{-1}>> C_*(K(n_1))\otimes\ldots \otimes C_*(K(n_m-j+1))\otimes C_*(K(j))\otimes\ldots \otimes C_*(K(n_k))\\
@VVV                @VV\id\otimes \circ_i \otimes\id V\\
F_*(\Gamma) @>\phi_{\Gamma}^{-1}>> C_*(K(n_1))\otimes\ldots\otimes C_*(K(n_m))\otimes\ldots\otimes C_*(K(n_k))\\ \end{CD}\]
Here the left vertical arrow is the map $p_{(\Gamma'\rightarrow\Gamma)}$.
\end{lem}

\proof Let $T_l\in C_*(K(n_l))$ $(1\leq l\leq k, l\neq m)$, $U\in C_*(K(n_m-j+1))$, and $V\in C_*(K(j))$. By definition of $\phi$ in Corollary~\ref{ch-map-cor}, we have
\[ \phi_{\Gamma'}^{-1}(\Gamma'(T_1,\ldots,U,V,\ldots,T_k)=(-1)^{|\Gamma'|} A_{\Gamma'}(-1)^{\clubsuit} T_1\otimes\ldots\otimes U\otimes V\otimes\ldots\otimes T_k,\]
where the sign is given by
\[\clubsuit= \sum_{1\leq p<q<m} n_p t_q + (n_1+\ldots+n_{m-1})t_U+(n_1+\ldots+n_m-j+1)t_V+\sum_{m<l} (n_1+\ldots+n_{l-1}+1) t_l.\]
Here, as before, $t_l$ is the number of internal vertices of $T_l$. 

Applying the morphism $(\phi_\Gamma)\circ (\id\otimes\circ_i\id)$, we obtain
\begin{align*}
[(\phi_\Gamma)&\circ (\id\otimes\circ_i\id)]\phi_{\Gamma'}^{-1}\big(\Gamma'(T_1,\ldots,U,V,\ldots,T_k)\big)=\\
&=\phi_\Gamma\big((-1)^{|\Gamma'|}A_{\Gamma'}(-1)^\clubsuit T_1\otimes\ldots\otimes U\circ_i V \otimes\ldots\otimes T_k\big)\\
&=(-1)^{|\Gamma'|}A_{\Gamma'}(-1)^\clubsuit (-1)^{|\Gamma|} A_\Gamma (-1)^{\sum_{1\leq p<q\leq m} n_p t_q} \Gamma(T_1,\ldots,U\circ_i V,\ldots, T_k)\\
&=(-1)^{|\Gamma'|}A_{\Gamma'}(-1)^\clubsuit (-1)^{|\Gamma|} A_\Gamma (-1)^{\sum_{1\leq p<q\leq m} n_p t_q}(-1)^\epsilon \Gamma'(T_1,\ldots,U, V,\ldots, T_k)
\end{align*}
where $\epsilon$ is as in Lemma~\ref{lema2}, and $t_m:=t_V+t_U+1$. We need to show that the sign in the last equation all cancel. For this we observe that
\[ \clubsuit+\sum_{1\leq p<q\leq m} n_p t_q= \sum_{l<m} n_l+ (n_m-j+1)t_V+\sum_{l>m} t_l \pmod 2,\]
and so
\[ \clubsuit+\sum_{1\leq p<q\leq m} n_p t_q+\epsilon= \sum_{l<m} n_l +i(j+1)+j(n_m-j+1) \pmod 2.\]
Finally, by the second part of Lemma~\ref{lema2} we have
 \[ A_\Gamma A_{\Gamma'}=(-1)^{\sum_{l<m} n_l +i(j+1)+j n_m}.\]
Hence we arrived at the desired conclusion that
\[(-1)^{|\Gamma'|}A_{\Gamma'}(-1)^\clubsuit (-1)^{|\Gamma|} A_\Gamma (-1)^{\sum_{1\leq p<q\leq m} n_p t_q}(-1)^\epsilon =1,\]
since $j(-j+1)=0\pmod 2$ and $|\Gamma'|=|\Gamma|+1$.\qed

\medskip
\noindent{\bf Proof of Proposition~\ref{compatible-prop}.} To simplify the notations, in this proof, we set
\[ C_{n_1,\ldots,n_k}:= C_*(K(n_1))\otimes\ldots\otimes C_*(K(n_k))\]
for a sequence of (positive) integers $n_1,\ldots,n_k$. Then the proposition follows from the following commutative diagram, where the middle square's commutativity is due to that $\Delta$ is an operad map of degree zero.
\[\begin{CD}
F_*(\Gamma') @>>> F_*(\Gamma)\\
@V\phi^{-1}_{\Gamma'}VV @VV\phi^{-1}_\Gamma V\\
C_{n_1,\ldots,n_m-j+1,j,\ldots,n_k}  @>\id\otimes\circ_i\otimes\id>> C_{n_1,\ldots,n_m,\ldots,n_k} \\
@V\Delta\otimes\ldots\otimes\Delta VV    @VV\Delta\otimes\ldots\otimes\Delta V        \\
C_{n_1,n_1,\ldots,n_m-j+1,n_m-j+1,j,j,\ldots,n_k,n_k}   @. C_{n_1,n_1,\ldots,n_m,n_m,\ldots,n_k,n_k}\\
 @V\tau' VV               @VV\tau V \\
(C_{n_1,\ldots,n_m-j+1,j,\ldots,n_k})^{\otimes 2} @>(\id\otimes\circ_i\otimes\id)^{\otimes 2}>> (C_{n_1,\ldots,n_m,\ldots,n_k})^{\otimes 2} \\
@V(\phi_{\Gamma'})^{\otimes 2} VV    @VV(\phi_\Gamma)^{\otimes 2} V\\
F_*(\Gamma')^{\otimes 2} @>>> F_*(\Gamma)^{\otimes 2}
\end{CD}\]
Here the top and bottom horizontal arrows are inclusions. The left vertical composition is by definition $\delta_{\Gamma'}$, while the right vertical composition is by definition $\delta_\Gamma$. Furthermore, the commutativity of the top and bottom squares follows from Lemma~\ref{inclusion-lem}. The commutativity of the middle square is due to Koszul sign property.\qed

\begin{cor}~\label{square}
The following diagram is commutative.
\[\begin{CD}
F_*(\Gamma) @> \delta_\Gamma >> F_*(\Gamma)\otimes F_*(\Gamma)\\
@V p_\Gamma VV          @VV p_\Gamma\otimes p_\Gamma V\\
\mG_* @>\delta >> \mG_*\otimes\mG_*.
\end{CD}\]
\end{cor}

\proof Let $\big(\Gamma_m\rightarrow \Gamma\big)$ be an element of $F_*(\Gamma)$. By factoring this morphism into several morphism each one contracting a single edge and iterated use of Proposition~\ref{compatible-prop}, we conclude there is a commutative diagram
\[\begin{CD}
F_*(\Gamma_m) @>>> F_*(\Gamma)\\
@V \delta_{\Gamma_m} VV             @V\delta_{\Gamma} VV\\
F_*(\Gamma_m)\otimes F_*(\Gamma_m) @>>> F_*(\Gamma)\otimes F_*(\Gamma)
\end{CD}\]
The element $\big(\Gamma_m\rightarrow \Gamma\big)\in F_*(\Gamma)$ lifts to an element $\big(\Gamma_m\stackrel{\id}{\rightarrow} \Gamma_m\big)\in F_*(\Gamma_m)$. By commutativity of the above diagram and naturality of $F_*$, it suffices to prove that
\[ \delta p_{\Gamma_m} \Big(\big(\Gamma_m\rightarrow \Gamma_m\big)\Big) = p_{\Gamma_m}\otimes p_{\Gamma_m} \delta_{\Gamma_m}\Big(\big(\Gamma_m\rightarrow \Gamma_m\big)\Big).\]
By definition of $\delta_{\Gamma_m}$, we have
\begin{align*}
\delta_{\Gamma_m}\Big(\big(\Gamma_m\rightarrow \Gamma_m\big)\Big)&=(-1)^{|\Gamma_m|}A_{\Gamma_m} \cdot (-1)^{\sum_{i<j}(n_i+t_{n_i}^{(2)})(n_j+t_{n_j}^{(1)})}\cdot (-1)^{\sum_{i<j}n_it_{n_j}^{(1)}+n_it_{n_j}^{(2)}} \cdot \\
&\big(\Gamma_m(c_{n_1}^{(1)},\ldots,c_{n_k}^{(1)})\rightarrow \Gamma_m\big)\otimes \big(\Gamma_m (c_{n_1}^{(2)},\ldots,c_{n_k}^{(2)})\rightarrow \Gamma_m\big).\end{align*}
We note that the middle sign $(-1)^{\sum_{i<j}(n_i+t_{n_i}^{(2)})(n_j+t_{n_j}^{(1)})}$ comes from the Koszul sign of the permutation $\tau$ in Definition~\ref{delta-gamma-def}. Using the fact that $t_{n_i}^{(1)}+t_{n_i}^{(2)}=n_i\pmod 2$, the above sign simplifies to
\begin{align*}
\;\;\;\;\;\;\;\;\;\;& (-1)^{|\Gamma_m|}A_{\Gamma_m} \cdot (-1)^{\sum_{i<j} t_{n_i}^{(1)}t_{n_j}^{(2)}+n_it_{n_j}^{(1)}+n_it_{n_j}^{(2)}}\\
 &=  (-1)^{|\Gamma_m|}A_{\Gamma_m} \cdot (-1)^{\sum_{i<j}n_it_{n_j}^{(1)}+(t_{n_i}^{(1)}+n_i)t_{n_j}^{(2)}}\\
 &=(-1)^{|\Gamma_m|}A_{\Gamma_m} \cdot (-1)^{\sum_{i<j}n_it_{n_j}^{(1)}+t_{n_i}^{(2)}t_{n_j}^{(2)}}
 \end{align*}
 Now applying the map $p_{\Gamma_m}\otimes p_{\Gamma_m}$ gives exactly $\delta (\Gamma_m)$, which equals $\delta(p_{\Gamma_m}(\Gamma_m))$\qed
 
\begin{cor}
The linear map $\delta: \mG_*\rightarrow \mG_*\otimes \mG_*$ is a morphism of chain complexes.
\end{cor}

\proof Let $\Gamma\in \mG_*$, by the previous corollary, we may lift $\Gamma$ to $(\Gamma\rightarrow\Gamma)\in F_*(\Gamma)$. Then we have
\begin{align*}
\partial \delta(\Gamma)&= \partial \delta p_\Gamma\big((\Gamma\rightarrow\Gamma)\big)\\
&=\partial p_\Gamma\otimes p_\Gamma \circ \delta_\Gamma\big((\Gamma\rightarrow\Gamma)\big)\\
&=p_\Gamma\otimes p_\Gamma \circ \delta_\Gamma\big(\partial (\Gamma\rightarrow\Gamma)\big)\\
&=\delta p_\Gamma\big(\partial (\Gamma\rightarrow\Gamma)\big)\\
&=\delta \partial \Gamma.
\end{align*}
Here we used the fact that $p_\Gamma$ and $\delta_\Gamma$ are chain maps.\qed

The following proposition proves that $\delta$, under the rational homotopy equivalence $\phi: C_*(\mathcal{F}at)\rightarrow \mG_*$, is homotopic to the Alexander-Whitney diagonal on $C_*(\mathcal{F}at)$.

\medskip
\begin{prop}~\label{diagonal-prop}
Let $\phi: C_*(\mathcal{F}at) \rightarrow \mG_*$ be the unique morphism carried by the forest carrier $F_*$. Denote by $\Delta_{AW}$ the Alexander-Whitney diagonal on $C_*(\mathcal{F}at)$. Then the following diagram is commutative up to homotopy.
\[\begin{CD}
C_*(\mathcal{F}at) @>\Delta_{AW}>> C_*(\mathcal{F}at)\otimes C_*(\mathcal{F}at)\\
@V\phi VV    @V\phi\otimes\phi VV\\
\mG_* @>\delta>> \mG_*\otimes \mG_*
\end{CD}\]
\end{prop}

\proof Let $\mathcal{F}in$ be the full subcategory of $\mathcal{F}at$ which contains exactly one object from every isomorphism class. In~\cite{Igu}, Igusa constructed a chain map
\[ \psi: \mG_* \rightarrow C_*(\mathcal{F}in)\]
which is rational inverse to the morphism $\phi$, after composing with the deformation retract $C_*(\mathcal{F}in)\hookrightarrow C_*(\mathcal{F}at)$. Thus to prove the proposition, it is enough to show that the two morphisms
\begin{align*}
C_*(\mathcal{F}in) &\stackrel{\Delta_{AW}}{\longrightarrow} C_*(\mathcal{F}in)\otimes C_*(\mathcal{F}in)\\
C_*(\mathcal{F}in) &\stackrel{\Delta^\dagger}{\longrightarrow} C_*(\mathcal{F}in)\otimes C_*(\mathcal{F}in)
\end{align*}
are homotopic, where the second morphism $\Delta^\dagger$ is the composition
\[ C_*(\mathcal{F}in) \hookrightarrow C_*(\mathcal{F}at) \stackrel{\phi}{\longrightarrow} \mG_*\stackrel{\delta}{\longrightarrow} \mG_*\otimes\mG_* \stackrel{\psi\otimes\psi}{\longrightarrow} C_*(\mathcal{F}in)\otimes C_*(\mathcal{F}in).\]
For this, we use the identity carrier from the category $\mathcal{F}in$ to the complex $C_*(\mathcal{F}in)$, as described in~\cite[Lemma 1.29]{Igu}. Recall the identity carrier $I$ is defined by assigning an object $\Gamma\in \mathcal{F}in$ to the chain complex $C_*(\mathcal{F}in/\Gamma)$ where $\mathcal{F}in/\Gamma$ is the over category of $\Gamma$. 
Moreover we have a morphism
\[ q_\Gamma: I(\Gamma)=C_*(\mathcal{F}in/\Gamma) \rightarrow C_*(\mathcal{F}in)\]
induced by the canonical functor $\mathcal{F}in/\Gamma \rightarrow \mathcal{F}in$ which ``forgets" being over $\Gamma$. Note that since $\mathcal{F}in/\Gamma$ admits a final object, the simplicial complex $C_*(\mathcal{F}in/\Gamma)$ is acyclic. Thus $I$ forms an acyclic carrier. We consider the tensor product carrier $I\otimes I$ defined by
\[ (I\otimes I) (\Gamma)= C_*(\mathcal{F}in/\Gamma)\otimes C_*(\mathcal{F}in/\Gamma),\]
which by K\"unneth is again acyclic.
To complete the proof of the proposition it is enough to show that both $\Delta_{AW}$ and $\Delta^\dagger$ are carried by $I\otimes I$, since any two morphisms carried by one acyclic carrier are homotopic.

We first show that the diagonal $\Delta_{AW}$ is carried by $I\otimes I$. Indeed for an $n$-simplex $(\Gamma_0\rightarrow\cdots\rightarrow\Gamma_n)$, define
\[ \widetilde{\Delta_{AW}}(\Gamma_0\rightarrow\cdots\rightarrow\Gamma_n) := \sum_{i=0}^n (\Gamma_0\rightarrow\cdots\rightarrow\Gamma_i)\otimes (\Gamma_i\rightarrow\cdots\rightarrow\Gamma_n)\in C_*(\mathcal{F}in/\Gamma_n)^{\otimes 2}.\]
Here we consider the sequence $(\Gamma_0\rightarrow\cdots\rightarrow\Gamma_i)$ in $C_*(\mathcal{F}in/\Gamma_n)$ via the composition $\Gamma_i\rightarrow \Gamma_{i+1}\cdots\rightarrow \Gamma_n$. It is immediate from the definition that $\widetilde{\Delta_{AW}}(\Gamma_0\rightarrow\cdots\rightarrow\Gamma_n)$ gives a required lifting of $\Delta_{AW}$.

Next we show that $\Delta^\dagger$ is also carried by $I\otimes I$. We claim the following diagram commutes
\[ \begin{CD}
F_*(\Gamma_n) @>\delta_{\Gamma_n}>> F_*(\Gamma_n)\otimes F_*(\Gamma_n) @>\widetilde{\psi}_{\Gamma_n}\otimes\widetilde{\psi}_{\Gamma_n} >> C_*(\mathcal{F}in/\Gamma_n)\otimes C_*(\mathcal{F}in/\Gamma_n)\\
@Vp_{\Gamma_n} VV  @VV p_{\Gamma_n}\otimes p_{\Gamma_n} V   @VV q_{\Gamma_n}\otimes q_{\Gamma_n} V\\
\mG_* @>\delta>> \mG_*\otimes\mG_* @>\psi\otimes\psi>> C_*(\mathcal{F}in)\otimes C_*(\mathcal{F}in).
\end{CD}\]
In fact each square commutes: Corollary~\ref{square} implies that the left square commutes, while the proof of~\cite[Lemma 1.29]{Igu} implies that the right square commutes.
Since $\phi$ is carried by $F_*$, for each $n$-simplex $(\Gamma_0\rightarrow\cdots\rightarrow\Gamma_n)$, there exists a lift
\[ \widetilde{\phi}(\Gamma_0\rightarrow\cdots\rightarrow\Gamma_n)\in F_*(\Gamma_n), \;\;\mbox{such that}\;\; p_{\Gamma_n}\widetilde{\phi}(\Gamma_0\rightarrow\cdots\rightarrow\Gamma_n)=\phi(\Gamma_0\rightarrow\cdots\rightarrow\Gamma_n).\]
Define $\widetilde{\Delta^\dagger}(\Gamma_0\rightarrow\cdots\rightarrow\Gamma_n):=(\widetilde{\psi}_{\Gamma_n}\otimes\widetilde{\psi}_{\Gamma_n})
\delta_{\Gamma_n}\widetilde{\phi}(\Gamma_0\rightarrow\cdots\rightarrow\Gamma_n)\in C_*(\mathcal{F}in/\Gamma_n)\otimes C_*(\mathcal{F}in/\Gamma_n)$. Then by the above commutative diagram we have
\begin{align*}
&(q_{\Gamma_n}\otimes q_{\Gamma_n})\widetilde{\Delta}(\Gamma_0\rightarrow\cdots\rightarrow\Gamma_n)\\
&=(q_{\Gamma_n}\otimes q_{\Gamma_n})(\widetilde{\psi}_{\Gamma_n}\otimes\widetilde{\psi}_{\Gamma_n})
\delta_{\Gamma_n}\widetilde{\phi}(\Gamma_0\rightarrow\cdots\rightarrow\Gamma_n)\\
&=(\psi\otimes\psi)\delta\phi(\Gamma_0\rightarrow\cdots\rightarrow\Gamma_n)\\
&=\Delta^\dagger(\Gamma_0\rightarrow\cdots\rightarrow\Gamma_n).
\end{align*}
The naturality of $\widetilde{\Delta^\dagger}$ follows from Proposition~\ref{compatible-prop}.\qed

Since the (dual of the) Alexander-Whitney diagonal defines the cup product on cohomology, this proposition immediately implies the following

\begin{cor}\label{cup-cor}
Given maps $c_1,c_2:\mG_*\lto \mQ$ representing cohomology classes $$[c_1],[c_2]\in H^*(\mG_*,\mQ)\simeq H^*\left(\coprod_{g,n}\mM_{g,n},\mQ\right),$$ we have
$[c_1]\cup[c_2]=[(c_1\otimes c_2)\circ \delta].$
\end{cor}

\subsection{Kontsevich classes} Let  $(A,\rho,\langle ,\rangle)$ be a finite dimensional cyclic $A_\infty$ algebra over $\mQ$. In this subsection, we shall assume that the inner product $\langle-,-\rangle$ on $A$ is even (and symmetric by our convention). The odd case will be dealt with in the next section. In~\cite{Kon}, Kontsevich constructed a cohomology class 
$$[c_A]\in H^*(\mG_*,\mQ)\cong H^*(\coprod_{g,n}\mM_{g,n},\mQ).$$ It is known this is a homotopy invariant of $A$, see~\cite{HamLaz}.

We first recall the definition of $c_A: \mG_*\rightarrow \mQ$, following~\cite[Section 2]{Igu}. In {\cite{Igu}}, Igusa wrote down the formula for $c_A$ explicitly by choosing a basis of $A$, and dealing carefully with signs involved. We shall use a more diagrammatic approach instead. The equivalence between the approaches will be clear from our construction. A notable feature of our definition is that we do not need to assume that $m_1=0$~\footnote{This is not a big generalization since we assume finite dimensionality. In this case, cyclic $A_\infty$ algebras with non-vanishing $m_1$ always admits a self-adjoint homotopy. Using the tree formula to transfer the $A_\infty$ structure, we obtain a homotopy equivalent cyclic $A_\infty$-algebra with $m_1=0$.}.  

\begin{defn}~\label{def-kont} Let $\Gamma\in \mG_*$ be a ribbon graph endowed with an orientation. We define a chain map $\widetilde{c_{A,\Gamma}}: F_*(\Gamma)\rightarrow \mQ$ as the composition of the maps in the following diagram
\[ \begin{CD}
F_*(\Gamma) @> \phi_\Gamma^{-1}>> C_*(K(n_1))\otimes\ldots\otimes C_*(K(n_k))\\
@.                  @VV\rho\otimes\ldots\otimes\rho V\\
          @.     \Hom(A^{\otimes n_1}, A)\otimes\ldots\otimes\Hom(A^{\otimes n_k},A) \\
@.           @VVD_k V\\
 @.            A^{\otimes n_1 +1} \otimes\ldots\otimes A^{\otimes n_k+1}\\
@.          @VV\eta V\\
@.      A^{\otimes 2}\otimes\ldots\otimes A^{\otimes 2}\\
@.         @VV\langle-,-\rangle\otimes\ldots\otimes\langle-,-\rangle V\\
@.       \mQ \otimes\ldots\otimes \mQ \cong \mQ.
\end{CD}\]
The vertical isomorphism $D_k$ in this diagram is defined as the tensor product of the compositions
\[\Hom(A^{\otimes n_l}, A) \xrightarrow{\theta^{-1}} A\otimes (A^\vee)^{\otimes n_l} \xrightarrow{\id\otimes D^{\otimes n_l}} A\otimes A^{\otimes n_l}=A^{\otimes n_l+1}, \ \textrm{for} \ l=1,\ldots,k\]
where $D: A^\vee\rightarrow A$ is the isomorphism induced by the non-degenerate pairing $\langle-,-\rangle$ and $\theta$ is the isomorphism defined by 
$$\theta(a \otimes \alpha_{n_l} \otimes \ldots \otimes \alpha_1)(x_1,\ldots,x_{n_l})=\alpha_1(x_1)\ldots\alpha_{n_l}(x_{n_l})a.$$
Secondly, the isomorphism $\eta$ is defined as follows. We fix an ordering of edges of $\Gamma$, as well as an orientation for each edge, and define $\eta$ to be the permutation associated to the permutation of half edges
\[ (e_{1,0},e_{1,n_1}\ldots,e_{1,1},\ldots,e_{k,0},e_{k,n_k},\ldots,e_{k,1}) \rightarrow ( h_1^+,h_1^-, \ldots, h_E^+,h_E^-)\]
where E is the number of edges of $\Gamma$ and $(-)^+$ and $(-)^-$ are the two half-edges of an edge.
\end{defn}

\begin{lem}\label{ind.class}
 The map $\widetilde{c_{A,\Gamma}}$ is independent of all the choices made, namely, orderings of vertices and edges of $\Gamma$, orientations of each edge and choices of roots $e_{i,0}$.
\end{lem}
\proof
The independence of the choices of ordering of vertices and roots $e_{i,0}$ is proved in the same way as Lemma \ref{lema.delta}, this time using the cyclicity of $\rho$.

Next note that the morphism $\langle-,-\rangle\otimes\ldots\otimes\langle-,-\rangle$ is invariant under permutations since $\langle-,-\rangle$ is even, therefore the composition $\widetilde{c_{A,\Gamma}}$ is independent of ordering of edges. Also, since $\langle -,- \rangle$ is symmetric, $\widetilde{c_{A,\Gamma}}$ does not depend on the choices of orientation of each edge.
\qed

\begin{prop}\label{c-tilde}
 Let $\Gamma$ and $\Gamma'$ be as in Proposition \ref{compatible-prop}. Then the following diagram is commutative.
\[\begin{CD}
F_*(\Gamma') @>p_{(\Gamma'\rightarrow\Gamma)}>> F_*(\Gamma)\\
@V \widetilde{c_{A,\Gamma'}} VV             @V \widetilde{c_{A,\Gamma}} VV\\
\mQ @= \mQ
\end{CD}\]
\end{prop}
\proof
The proof is similar to the proof of Proposition \ref{compatible-prop}. It follows from the commutativity of the following diagram
\[\begin{CD}
F_*(\Gamma') @>>> F_*(\Gamma)\\
@V\phi^{-1}_{\Gamma'}VV @VV\phi^{-1}_\Gamma V\\
C_{n_1,\ldots,n_m-j+1,j,\ldots,n_k}  @>\id\otimes\circ_i\otimes\id>> C_{n_1,\ldots,n_m,\ldots,n_k} \\
@V\rho^{\otimes k+1} VV    @VV\rho^{\otimes k} V        \\
\End^{n_1}\otimes\ldots \End^{n_m-j+1}\otimes\End^{j}\ldots \End^{n_k}   @>\id\otimes\circ_i\otimes\id>> \End^{n_1}\otimes\ldots \otimes \End^{n_k}\\
 @VD_{k+1} VV               @VV D_k V \\
A^{\otimes n_1+1}\otimes\ldots A^{\otimes n_m-j+2}\otimes A^{\otimes j+1}\ldots A^{\otimes n_k+1} @. A^{\otimes n_1+1}\otimes\ldots \otimes A^{\otimes n_k+1} \\
@V\eta' VV               @VV\eta V \\
A^{\otimes 2}\otimes\ldots\otimes A^{\otimes 2} @>\langle -,- \rangle \otimes \id>> A^{\otimes 2}\otimes\ldots\otimes A^{\otimes 2}\\
@V\langle -,- \rangle^{\otimes E+1} VV    @VV\langle -,- \rangle^{\otimes E} V\\
\mathbb{Q} @= \mathbb{Q}
\end{CD}\]
where we use the notation $\End^n:=\Hom(A^{\otimes n},A)$. The first square commutes by Lemma \ref{inclusion-lem}, the second commutes because $\rho$ is an operad map and the last square commutes because the inner product has degree zero. Finally commutativity of the third square is a straightforward computation that we omit.
\qed

\begin{defn}
Let $A$ be a finite dimensional cyclic $A_\infty$ algebra whose inner product is even. Define its Kontsevich class $c_A: \mG_*\rightarrow \mQ$ by
\[ c_A(\Gamma):=\widetilde{c_{A,\Gamma}}\big((\Gamma\stackrel{\id}{\rightarrow}\Gamma)\big).\]
\end{defn}
One can check that this definition agrees with the one given in \cite{Igu}.

\begin{lem}
The linear map $c_A: \mG_*\rightarrow \mQ$ is a chain map, i.e. $c_A$ is a cocycle of $\mG_*$.
\end{lem}

\proof Let $\partial\Gamma=\sum_{\Gamma'}\Gamma'$ in $\mG_*$. By applying the previous proposition, we have
\begin{align*}
c_A(\partial\Gamma)&=\sum_{\Gamma'} c_A(\Gamma')\\
&= \sum_{\Gamma'} \widetilde{c_{A,\Gamma'}}\big((\Gamma'\stackrel{\id}{\rightarrow}\Gamma')\big)= \sum_{\Gamma'} \widetilde{c_{A,\Gamma}}\big((\Gamma'\rightarrow\Gamma)\big)\\
&=\widetilde{c_{A,\Gamma}}\big(\sum_{\Gamma'} (\Gamma'\rightarrow\Gamma)\big)=\widetilde{c_{A,\Gamma}}\big(\partial (\Gamma\stackrel{\id}{\rightarrow}\Gamma)\big)\\
&=\partial\widetilde{c_{A,\Gamma}}\big((\Gamma\stackrel{\id}{\rightarrow}\Gamma)\big)=0.
\end{align*}
The last line follows from the fact that $\widetilde{c_{A,\Gamma}}$ is a chain map.\qed

\begin{thm}~\label{tensor-formula-thm}
Let $A$ and $B$ be two finite dimensional cyclic $A_\infty$ algebras. Assume that the inner products on both $A$ and $B$ are even. Fix a diagonal $\Delta$ of $\Aii$ to define the cyclic tensor product $A\otimes B$ and the morphism $\delta: \mG_*\rightarrow\mG_*\otimes \mG_*$. Then we have
\[ c_{A\otimes B}= (c_A\otimes c_B)\circ \delta.\]
In view of Corollary~\ref{cup-cor}, this implies $[c_{A\otimes B}]= [c_A]\cup [c_B]$.
\end{thm}

\proof As before we denote
\[ C_{n_1,\ldots,n_k}:= C_*(K(n_1))\otimes\ldots\otimes C_*(K(n_k)),\]
and $\End^n_V:=\Hom(V^{\otimes n},V)$ for a vector space $V$. 

First observe that for an oriented ribbon graph $\Gamma\in \mG_*$, Proposition \ref{c-tilde} and Corollary~\ref{square} imply that
\[ (c_A\otimes c_B)\circ\delta(\Gamma)=(\widetilde{c_{A,\Gamma}}\otimes\widetilde{c_{B,\Gamma}})\circ\delta_\Gamma \big((\Gamma\stackrel{\id}{\rightarrow}\Gamma)\big).\]
Thus it suffices to prove that
\[ \widetilde{c_{A\otimes B, \Gamma}}=(\widetilde{c_{A,\Gamma}}\otimes \widetilde{c_{B,\Gamma}})\circ (\delta_\Gamma)\]
for a fixed $\Gamma\in \mG_*$. The proof of this equality follows from commutativity of the following diagram.

\medskip
\[\begin{CD}
F_*(\Gamma) @= F_*(\Gamma)\\
@V\phi_\Gamma^{-1} VV @VV\phi_\Gamma^{-1} V\\
C_{n_1,\ldots,n_k}  @=  C_{n_1,\ldots,n_k}\\
           @V\Delta\otimes\ldots\otimes\Delta VV         @VV\Delta\otimes\ldots\otimes\Delta V\\
     C_{n_1,n_1,\ldots,n_k,n_k} @=  C_{n_1,n_1,\ldots,n_k,n_k}\\
     @| @VV\tau V\\
     C_{n_1,n_1,\ldots,n_k,n_k} @>\tau>> C_{n_1,\ldots,n_k}\otimes C_{n_1,\ldots,n_k}\\
 @V(\rho_A\otimes\rho_B)\otimes\ldots\otimes(\rho_A\otimes\rho_B)VV   @VV(\rho_A\ldots\rho_A)\otimes(\rho_B\ldots\rho_B) V\\
\End^{n_1}_A\otimes\End^{n_1}_B\ldots\End^{n_k}_A\otimes\End^{n_k}_B @>\tau>> \End^{n_1}_A\ldots\End^{n_k}_A\otimes\End^{n_1}_B\ldots\End^{n_k}_B\\
  @Vi\otimes\ldots\otimes i VV         @|\\
\End^{n_1}_{A\otimes B}\ldots\End^{n_k}_{A\otimes B}  @.    \End^{n_1}_A\ldots\End^{n_k}_A\otimes\End^{n_1}_B\ldots\End^{n_k}_B\\
 @V\cong VV            @VV\cong V \\
(A\otimes B)^{n_1+1}\ldots (A\otimes B)^{n_k+1} @>\tau >> A^{n_1+1}\ldots A^{n_k+1}\otimes B^{n_1+1}\ldots B^{n_k+1}\\
@V\eta_{A\otimes B} VV     @VV\eta_A\otimes \eta_B V\\
(A\otimes B)^{\otimes 2E}   @>\tau>> A^{\otimes 2E} \otimes B^{\otimes 2E}\\
@V(\langle-,-\rangle_{A\otimes B})^{\otimes E} VV   @VV(\langle-,-\rangle_A)^{\otimes E}\otimes (\langle-,-\rangle_B)^{\otimes E} V\\
\mQ \cong \mQ^{\otimes E} @>\cong >> \mQ^{\otimes E}\otimes \mQ^{\otimes E}=\mQ.\\
\end{CD}\]

\medskip
\noindent Here the horizontal $\tau$'s are the appropriate permutations. Observe that the left vertical composition gives $\widetilde{c_{A\otimes B,\Gamma}}$, while the right vertical composition gives $(\widetilde{c_{A,\Gamma}}\otimes \widetilde{c_{B,\Gamma}})\circ \delta_\Gamma$. The top three squares obviously commute and the fourth ones commutes by Koszul sign convention and the fact that $\rho$ is a map of degree zero. The bottom square commutes by the definition of $\langle-,-\rangle_{A\otimes B}$. The commutativity of the other squares simply follows from the Koszul sign convention and the definitions.  \qed

\section{Twisted cases}\label{sec:twisted} 
\subsection{Twisted ribbon graph complex} In the case when $A$ is a cyclic $A_\infty$-algebra with an \emph{odd} cyclic inner product, one can still define a Kontsevich class $[c_A]$ in the cohomology of $\coprod_{g,n} \mM_{g,n}$, but now with coefficients in a local system $\LL$. The fiber of $\LL$ over a Riemann surface $\Sigma$ is the determinant (or top exterior power) of $H_1(\Sigma)$.  In our model for $\coprod_{g,n} \mM_{g,n}$, the nerve of the category $\mathcal{F}at$, the fiber of $\LL$, over a ribbon graph $\Gamma$, is
\[ \LL|_{\Gamma}:=\det H_1(\Gamma),\]
the determinant of the first homology group of $\Gamma$ (seen as a $1$-dimensional cell complex). Recall that morphisms in $\mathcal{F}at$ are defined as contractions of edges, hence they induce isomorphisms on homology and so define maps between the corresponding fibers of $\LL$. This data determines a local system in $C_*\mathcal{F}at$. Also note that since the group $\det H_1(\Gamma)$ is generated by an integral basis, it follows that $\LL^{\otimes 2}$ is canonically trivial. This implies that $\LL$ and its dual $\LL^\vee$ are canonically identified. Therefore we have
\[H^*\big(\coprod_{g,n} \mM_{g,n}, \LL\big)\cong H^*(C_*(\mathcal{F}at, \LL)^\vee),\
\]
where $C_*(\mathcal{F}at, \LL)^\vee$ is the dual of the simplicial complex of the nerve of $\mathcal{F}at$ with coefficients in $\LL$. We refer the reader to \cite{Spa} for basics on local systems.

In terms of the graph complex, twisting by $\LL$ amounts to a changing the definition of orientation for a graph. A good reference for this is \cite{LazVor}. Indeed, given a ribbon graph $\Gamma$ a twisted orientation $\mu$ on $\Gamma$ is defined to be a (integral) generator of
\[ \det E(\Gamma),\]
that is, an orientation on the vector space spanned by the set of edges of $\Gamma$. We will refer to such an orientation as a \emph{twisted orientation}.  
\begin{defn}
The twisted ribbon graph complex $\mG_*(\LL)$ is defined as the span of isomorphism classes of pairs $(\Gamma,\mu)$ where $\mu$ is a twisted orientation on $\Gamma$, modulo the relation $-(\Gamma,\mu)=(\Gamma,-\mu)$. The degree of $(\Gamma,\mu)$ is defined as $|\Gamma|=\sum_{v\in V(\Gamma)} \val(v)-3$. The differential is again given by summing  over expansion of vertices:
\[\partial(\Gamma,\mu):=\sum (\Gamma',\mu')\] 
The two orientations are related as $\mu'=e\wedge \mu$, where $e$ is the unique new edge of $\Gamma'$.
\end{defn}

We note that, since $\Gamma$ is a connected graph, there is a canonical identification
\begin{equation}~\label{eq:or}
 \det E(\Gamma) \cong \det(V(\Gamma)\oplus H(\Gamma))\otimes\det H_1(\Gamma).
\end{equation}
We refer the reader to~\cite{CV} for a proof of this fact. In the following we shall freely use this identification.

To be able to write down the twisted diagonals and odd Kontsevich classes, we first define a twisted version of Igusa's forest carrier on the category $\mathcal{F}at$. This will also allow us to compare $C_*(\mathcal{F}at,\LL)$ with $\mG(\LL)$.

For a ribbon graph $\Gamma$, we define the twisted forest carrier $F_*(\Gamma,\LL)$ as the span of isomorphism classes of ribbon graph morphisms $[\Gamma'\rightarrow \Gamma]$ together with $\mu'$ a twisted orientation on $\Gamma'$. Two such classes $\Gamma'\rightarrow \Gamma$ and $\Gamma''\rightarrow \Gamma$ are isomorphic if there is an twisted orientation preserving isomorphism $\Gamma'\rightarrow\Gamma''$ compatible with the maps to $\Gamma$. Observe that
\begin{itemize}
\item[(a)] The complex $F_*(\Gamma,\LL)$ does not have a natural augmentation, but a choice of a generator of $\xi$ of $\det H_1(\Gamma)$ induces one. Given $\Gamma'\rightarrow \Gamma \in F_0(\Gamma,\LL)$, $\xi$ determines a generator $\xi'$ of $\det H_1(\Gamma')$ since $\Gamma'\rightarrow \Gamma$ is a homotopy equivalence. Using (\ref{eq:or}) we can compare $\mu'$ with $c(\Gamma')\otimes\xi'$ and define the augmentation in a similar fashion to the untwisted case.
\item[(b)] There is a chain map $p^\LL_\Gamma: F_*(\Gamma,\LL)\rightarrow \mG(\LL)$, which forgets the base graph $\Gamma$.
\item[(c)] Given a morphism $\Gamma_1\rightarrow \Gamma_2$, there is a chain map $F_*(\Gamma_1,\LL)\rightarrow F_*(\Gamma_2,\LL)$ by post-composing with the map $\Gamma_1\rightarrow \Gamma_2$.
\end{itemize}

Underlying all the constructions in the previous section was the chain isomorphism $\phi_\Gamma$ defined in Corollary \ref{ch-map-cor}. We will now define a twisted version of this map. 

\begin{defn}
Let $\Gamma$ be a ribbon graph, choose an ordering of its vertices and starting half-edges at each vertex as in Construction~\ref{delta-constr}. Moreover choose a generator $\xi\in \det H_1(\Gamma)$.  We define
\[\psi_\Gamma: C_*(K(n_1))\otimes\ldots\otimes C_*(K(n_k)) \rightarrow F_*(\Gamma,\LL)\] 
by the formula
\[ \psi_\Gamma(T_1,\ldots,T_k):=(-1)^{|\Gamma|(\chi(\Gamma)+1)+\sum_{i<j}n_in_j} (-1)^{\sum_{i<j} n_i t_j } [\Gamma(T_1,\ldots,T_k)\rightarrow \Gamma].\]
Here $\chi(\Gamma)$ is the Euler characteristic of $\Gamma$ and the ribbon graph $\Gamma(T_1,\ldots,T_k)$ is as in Construction~\ref{delta-constr}, except for the orientations which are defined as follows,
\begin{align*}
&\mu_{\Gamma(T_1,\ldots,T_k)}= \sigma_1 \wedge \ldots \wedge \sigma_k \wedge \mu_\Gamma,\\
\mu_\Gamma&= \big(v_1\wedge e_{10}\ldots e_{1n_1}\wedge v_2\wedge e_{20}\ldots e_{2n_2}\ldots v_k\wedge e_{k0}\ldots e_{kn_k}\big)\otimes \xi,
\end{align*}
here $\sigma_i$ is the orientation of $T_i$ and we use the isomorphism~(\ref{eq:or}).
\end{defn}

One can easily show that $\psi_\Gamma$ is a chain isomorphism just like in Corollary \ref{ch-map-cor}. In particular $C_*(\mathcal{F}at,\LL)$ is acyclic and therefore it carries a chain map. As these complexes are not augmented, we have to change the definition of carrier as follows. Instead of requiring the maps to respect augmentations we require the following: given $(\Gamma_0, \xi)$ a degree zero, integral generator of $C_0(\mathcal{F}at, \LL)$, we require that $\widetilde{\phi^\LL}(\Gamma_0,\xi) \in F_*(\Gamma_0,\LL)$ satisfies $\epsilon_\xi(\widetilde{\phi^\LL}(\Gamma_0,\xi))=1$, where $\epsilon_\xi$ is the augmentation induced by $\xi$. As $\epsilon_{-\xi}=-\epsilon_{\xi}$, this is well defined. With this modification we have the following

\begin{prop}
The twisted forest carrier $F_*(\Gamma,\LL)$ induces a unique (up to homotopy) chain map
\[ \phi^\LL: C_*(\mathcal{F}at,\LL) \rightarrow \mG_*(\LL),\]
from the simplicial complex of the nerve of $\mathcal{F}at$, with coefficients in $\LL$, to the twisted ribbon graph complex. Moreover this map is a homotopy equivalence. Hence we get a homotopy equivalence $$C_*\big(\coprod_{g,n} \mM_{g,n}, \LL\big)\cong  \mG_*(\LL). $$
\end{prop}
This proposition can be proved by following Igusa's proof of the same statement in the untwisted case, in \cite{Igu}. As all the constructions, including the definition of $\psi^\LL$ the homotopy inverse of $\phi^\LL$ are entirely analogous to the one in \cite{Igu}, simply replacing orientations of ribbon graphs by twisted ones, we omit the proof.

\subsection{Odd Kontsevich class} In this subsection we define the Kontsevich class of an odd cyclic $A_\infty$-algebra. We begin with the following
\begin{defn}\label{def-kont-odd}
Assume that $A$ is a cyclic $A_\infty$-algebra endowed with an odd inner product. Let $\Gamma$ be a ribbon graph, and let $\xi$ be a generator of $\det H_1(\Gamma)$. Define $\widetilde{c_{A,\Gamma}}: F_*(\Gamma,\LL)\rightarrow \mQ$ as the composition
\[\begin{CD}
F_*(\Gamma,\LL) @> \psi_\Gamma^{-1}>> C_*(K(n_1))\otimes\ldots\otimes C_*(K(n_k))\\
@.                  @VV\rho\otimes\ldots\otimes\rho V\\
          @.     \Hom(A^{\otimes n_1}, A)\otimes\ldots\otimes\Hom(A^{\otimes n_k},A) \\
@.           @VV D_k V\\
 @.            A^{\otimes n_1 +1} \otimes\ldots\otimes A^{\otimes n_k+1}\\
@.          @VV\eta V\\
@.      A^{\otimes 2}\otimes\ldots\otimes A^{\otimes 2}\\
@.         @VV\langle-,-\rangle\otimes\ldots\otimes\langle-,-\rangle V\\
@.       \mQ \otimes\ldots\otimes \mQ \cong \mQ.
\end{CD}\]
Here all maps are the same as in Definition~\ref{def-kont}, except that when defining the permutation $\eta$ of half edges
\[ (e_{1,0},e_{1,n_1}\ldots,e_{1,1},\ldots,e_{k,0},e_{k,n_k},\ldots,e_{k,1}) \rightarrow ( h_1^+,h_1^-, \ldots, h_E^+,h_E^-),\]
we need to choose an ordering of the edges of $\Gamma$ that is compatible with $\mu_\Gamma$.
\end{defn}

The definition of $c_{A,\Gamma}$ is independent of all choices made. Indeed, independence on the ordering of the vertices can be proved as in the untwisted case and cyclicity of $\rho$ implies the independence on choices of starting edges $v_{i,0}$ (see Lemma \ref{ind.class}). Independence on orientation of each edge $(h_i^+,h_i^-)$ is due to the fact that the pairing is symmetric. Finally if we change $\xi$ to $-\xi$, the morphism $\psi_\Gamma$ changes by $-1$, which cancels the sign change in $\eta$ caused by interchanging two consecutive $\langle-,-\rangle$ due to the oddness of the inner product.

\medskip
Using $\widetilde{c_{A,\Gamma}}$, we define the Kontsevich class $c_A: \mG_*(\LL)\rightarrow \mQ$ by the formula
\begin{equation}\label{oddkontsevich}
c_A(\Gamma):= \widetilde{c_{A,\Gamma}}([\Gamma\stackrel{\id}{\rightarrow}\Gamma]),
\end{equation} 
where we take $\mu_\Gamma$ as the twisted orientation in $\Gamma$.

One can show that Lemma~\ref{inclusion-lem} in the twisted case remains valid.
This implies that the Kontsevich class $c_A: \mG_*(\LL) \rightarrow \mQ$ is a chain map. 

As in the even case~\cite{Igu}, one can write down an explicit formula of the Kontsevich class $c_A$ in terms of a summation over states by choosing a basis of the $A_\infty$-algebra $A$. Let $x_1,\cdots,x_d$ be a basis of $A$ and denote by $x_1^\vee,\cdots,x_d^\vee$, the dual basis. A state of a ribbon graph is a map of sets
\[ s: H(\Gamma) \rightarrow \left\{x_1,\cdots,x_n\right\},\]
from the set of half edges of the graph to the set of basis elements. Given a state, we denote the image of a half edge $e_{ij}$ by $x_{ij}$. Then the formula of $c_A$ is given by
\[ c_A(\Gamma)=\sum_{\mbox{states of }\Gamma} \epsilon\cdot B \prod_{\mbox{vertices}}\langle m_{n_i}(x_{i1},\cdots,x_{in_i}),x_{i0}\rangle\prod_{\mbox{edges}}\langle D((x^+)^\vee), D((x^-)^\vee)\rangle.\]
The signs $\epsilon$ and $B$ are given 
\begin{align*}
\epsilon &:=\prod_{\mbox{vertices}} (-1)^{n_i|x_{i0}|+(n_i-1)|x_{in_i}|+\cdots+1\cdot|x_{i2}|}\\
B&:= \mbox{the Koszul sign from the permutation $\eta$ used in Definition~\ref{def-kont-odd},}\\
&\mbox{\;\;\;\; with each half edge $e_{ij}$ labled by $D(x_{ij}^\vee)$.}
\end{align*}
Tracing back through the definitions one can check that this agrees with (\ref{oddkontsevich}). We will just point out that $c_A(\Gamma)$ vanishes unless $\Gamma$ has an even number of edges. This is because otherwise there are no states with non-trivial contributions, since the inner product is odd. This implies that the sign $(-1)^{|\Gamma|(\chi(\Gamma)+1)}$ in the definition of $\psi_\Gamma$, equals one since in this case $|\Gamma| = |\chi(\Gamma)| \pmod 2$.

\subsection{Twisted diagonals}
As we have seen, in the untwisted case, the proof of the tensor product formula for Kontsevich classes follows easily after the construction of a diagonal map $\delta$ on the graph complex. In the twisted situation, there are two cases to consider:
\begin{itemize}
\item[(1)] both $\langle-,-\rangle_A$ and $\langle-,-\rangle_B$ are odd;
\item[(2)] the inner product $\langle-,-\rangle_A$ is even, and $\langle-,-\rangle_B$ is odd (the case when $\langle-,-\rangle_A$ is odd and $\langle-,-\rangle_B$ is even is symmetric).
\end{itemize}
In the following, for a fixed graph $\Gamma$, as in Definition~\ref{delta-gamma-def}, we first define $\delta_\Gamma$ in each case separately.

\medskip
\noindent {\bf Case (1):}  We define $\delta_\Gamma: F_*(\Gamma,\LL) \rightarrow F_*(\Gamma,\LL)\otimes F_*(\Gamma,\LL)$ as the composition
\[\begin{CD}
F_*(\Gamma) @>\phi_\Gamma^{-1} >>     C_{n_1,\ldots,n_k} \\
        @.                                      @VV\Delta\otimes\ldots\otimes\Delta V\\
         @. C_{n_1,n_1,\ldots,n_k,n_k}\\
          @.                                       @VV\tau V\\
         F_*(\Gamma,\LL)\otimes F_*(\Gamma,\LL) @<\psi_\Gamma\otimes\psi_\Gamma<<   C_{n_1,\ldots,n_k} \otimes C_{n_1,\ldots,n_k}\\
\end{CD}\]


\medskip
\noindent {\bf Case (2):} We define $\delta_\Gamma: F_*(\Gamma,\LL) \rightarrow F_*(\Gamma)\otimes F_*(\Gamma,\LL)$ as  the composition
\[\begin{CD}
F_*(\Gamma,\LL) @>\psi_\Gamma^{-1} >>     C_{n_1,\ldots,n_k} \\
        @.                                      @VV\Delta\otimes\ldots\otimes\Delta V\\
         @. C_{n_1,n_1,\ldots,n_k,n_k}\\
          @.                                       @VV\tau V\\
         F_*(\Gamma)\otimes F_*(\Gamma,\LL) @<\phi_\Gamma\otimes\psi_\Gamma<<   C_{n_1,\ldots,n_k} \otimes C_{n_1,\ldots,n_k}\\
\end{CD}\]
In both cases the definition of $\delta_\Gamma$ does not depend on the choices involved. For ordering of vertices and starting half edges, the proof is the same as in the untwisted case. The map $\psi_\Gamma$ also depends on the choice of $\xi\in \det H_1(\Gamma)$, but note that the map $\psi_\Gamma$ (or its inverse) appears twice in both case, so this dependence cancels out.

Similar to the untwisted situation (Proposition~\ref{compatible-prop}), one can show that in the two twisted cases, $\delta_\Gamma$ are compatible with the maps $F_*(\Gamma',\LL) \rightarrow F_*(\Gamma,\LL)$ induced by a morphism $\Gamma'\rightarrow \Gamma$. This enables the definition of the twisted diagonal morphisms on the corresponding graph complexes.

\medskip
\noindent {\bf Case (1):} We define
\[ \delta: \mG_*(\LL)\rightarrow \mG_*\otimes\mG_*(\LL),\]
to be the unique morphism which makes the following diagram commutative (for all ribbon graph $\Gamma$)
\[\begin{CD}
F_*(\Gamma,\LL) @> \delta_\Gamma >> F_*(\Gamma)\otimes F_*(\Gamma,\LL)\\
@V p_\Gamma^\LL VV          @VV p_\Gamma\otimes p^\LL_\Gamma V\\
\mG_*(\LL) @>\delta >> \mG_*\otimes\mG_*(\LL).
\end{CD}\]
Or explicitly, for a ribbon graph $\Gamma$ endowed with a twisted orientation, we set
\[ \delta(\Gamma):=p_\Gamma\otimes p^\LL_\Gamma \big( \delta_\Gamma ([\Gamma\stackrel{id}{\rightarrow} \Gamma])\big).\]
This is well-defined by the twisted version of Proposition~\ref{compatible-prop} mentioned in the previous paragraph. 

\medskip
\noindent {\bf Case (2):} Similarly, we can define
\[ \delta: \mG_*\rightarrow \mG_*(\LL)\otimes\mG_*(\LL),\]
to make the diagram 
\[\begin{CD}
F_*(\Gamma) @> \delta_\Gamma >> F_*(\Gamma,\LL)\otimes F_*(\Gamma,\LL)\\
@V p_\Gamma VV          @VV p_\Gamma^\LL\otimes p^\LL_\Gamma V\\
\mG_* @>\delta >> \mG_*(\LL)\otimes\mG_*(\LL).
\end{CD}\]
commutative for all ribbon graph $\Gamma$.

In both cases we can prove the analogue of Proposition \ref{diagonal-prop}. That is, in Case (1) $\delta$ is homotopic, via $\phi^\LL$, to the map induced by the standard Alexander-Whitney diagonal
\[\Delta_{AW}:C_*(\mathcal{F}at,\LL)\rightarrow C_*(\mathcal{F}at)\otimes C_*(\mathcal{F}at,\LL),
\]
and similarly in Case (2). Therefore we can use $\delta$ to compute the cup product in cohomology.
Finally the proof of the tensor product formula is almost identical to that of Theorem~\ref{tensor-formula-thm}. We summarize these results in the following

\begin{thm}
Let $A$ and $B$ be two finite dimensional cyclic $A_\infty$-algebras whose inner product can be either even or odd. Then we have a chain level identity:
\[ c_{A\otimes B}= (c_A\otimes c_B)\circ \delta.\]
\end{thm}
This completes the proof of Theorem~\ref{main2-thm}.

\end{document}